\documentclass{article}
\usepackage{latexsym}
\usepackage{graphicx}
\renewcommand{\theequation}{\thesection.\arabic{equation}}
\newtheorem{theorem}{Theorem}
\newtheorem{lemma}{Lemma}
\begin{document}

\title{Maps close to identity and universal maps in the Newhouse domain}
\author{Dmitry Turaev\\
\small Imperial College, London, SW7 2AZ, UK\\ \small dturaev@imperial.ac.uk}
\maketitle

\begin{abstract}
Given an $n$-dimensional $C^r$-diffeomorphism $g$,
its renormalized iteration is an iteration of $g$, restricted to a certain
$n$-dimensional ball and taken in some $C^r$-coordinates in which the ball acquires
radius $1$. We show that for any $r\geq 1$ the renormalized iterations
of $C^r$-close to identity maps of an $n$-dimensional unit ball $B^n$ ($n\geq 2$)
form a residual set among all orientation-preserving $C^r$-diffeomorphisms $B^n\rightarrow R^n$.
In other words, any generic $n$-dimensional dynamical phenomenon
can be obtained by iterations of $C^r$-close to identity maps, with the same
dimension of the phase space. As an application, we show that any $C^r$-generic two-dimensional
map which belongs to the Newhouse domain (i.e., it has a wild hyperbolic set, so it is
not uniformly-hyperbolic, nor uniformly partially-hyperbolic) and which neither contracts, nor expands
areas, is $C^r$-universal in the sense that its iterations, after an appropriate coordinate transformation,
$C^r$-approximate every orientation-preserving two-dimensional diffeomorphism arbitrarily well.
In particular, every such universal map has an infinite set of coexisting hyperbolic attractors and
repellers.
\end{abstract}

\section{Ruelle-Takens problem and universal maps}
\label{rtu}\setcounter{equation}{0}
A long standing open problem in the theory of dynamical systems
is to describe which kind of dynamical phenomena can be expected
in close to identity maps. It started
with a much celebrated paper \cite{RT} where it was shown that any
$n$-dimensional dynamics can be implemented by a $C^n$-small perturbation
of the identity map of an $n$-dimensional torus. The paper seized a lot
of attention by physicists, because it proposed a new view on the
onset of hydrodynamical turbulence; at the same time it caused a
lot of criticism. One of the reasons for the critique was that the $C^n$-small
perturbations constructed in \cite{RT} were not small in
$C^{n+1}$, which is quite unphysical. The controversy was
resolved in \cite{NRT} where it was shown that for any $r$, given any
$C^r$-diffeomorphism $F$ of a closed $n$-dimensional ball, one can find
a $C^r$-close to identity map $g$ of the $(n+1)$-dimensional
closed unit ball $B^{n+1}$ such that the diffeomorphism
$F$ coincides with some iteration of the map $g$ restricted to some $n$-dimensional
invariant manifold. Thus, the restriction on smoothness
of perturbations was removed by sacrificing one dimension of the
phase space; anyway, other scenarios of the transitions
to turbulence had already been known.

From the purely mathematical point of view, the question still
remained unsolved: can an arbitrary $n$-dimensional dynamics be
obtained by iterations of a $C^r$-close to identity map of $B^n$,
i.e. in the same dimension of the phase space?
The difficulty is that the straightforward construction proposed in \cite{RT}
does not work for high $r$ in principle. Indeed, given an orientation-preserving
diffeomorphism $F:B^n\rightarrow R^n$, one can imbed it into a continuous family
${\cal F}_t$ of the diffeomorphisms such that ${\cal F}_1=F$ and ${\cal F}_0=id$.
Then, given any $N$, the map $F$ can be represented
as a superposition of $N$ maps
\begin{equation}\label{ruta}
F=F_{N}\circ \dots \circ F_1, \;\; \mbox{ where } \;\; F_s={\cal F}_{_{s/N}}\circ
{\cal F}_{_{(s-1)/N}}^{-1},
\end{equation}
that are $O(1/N)$-close to identity.
One can then choose $N$ pairwise disjoint small balls $D_s\in B^n$ of radius $\rho\sim N^{-1/(n-1)}$
and define a map $\phi: B^n\rightarrow B^n$ such that
$\phi(x)|_{x\in D_s}\equiv x_{s+1}+\rho F_s(\frac{x-x_s}{\rho})$
where $x_s$ is the center of $D_s$; the positions of the centers are chosen in such a way (Fig.1)
that the distances $\|x_{s+1}-x_s\|$ are uniformly close to zero for all $s$, hence
the map $\phi$ is $C^0$-close to identity. By construction, $\phi^{N}|_{D_0}$
is linearly conjugate to $F$, i.e.
the dynamics of $\phi^{N}|_{D_0}$ coincides with the dynamics of $F$. However, the derivatives of
$\phi$ of order $k$ behave as $N^{-1}\rho^{1-k}\sim N^{\frac{k-n}{n-1}}$, i.e. at
$k\geq n$ they do not, in general,
tend to zero as $N\rightarrow+\infty$. Thus, an arbitrary $n$-dimensional dynamics can be
implemented by iterations of $C^{n-1}$-close to identity maps of $B^n$, but the construction gives
no clue of whether the same can be said about the $C^n$-close to identity maps.

\begin{figure}
\centerline{\includegraphics[scale=0.4]{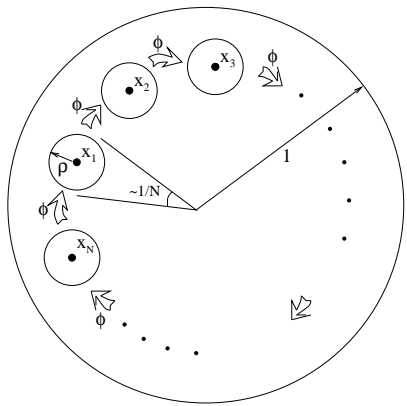}}
\caption{An illustration to Ruelle-Takens construction.}
\end{figure}

One could try to position the regions $D_s$ differently,
or make their radii vary, or change their shape. This, however, hardly can lead to an essential increase
in the maximal order of the derivatives (of $\phi$) which tend to zero as $N\rightarrow+\infty$.
The reason lies in a well-known fact from the averaging theory that the $O(\delta)$-close to identity map
$$\bar x= x+\delta f(x)$$
approximates a time shift of a certain autonomous flow with the accuracy
$O(\delta^m)$ for an arbitrarily large $m$ (if $f\in C^\infty$). Hence, the number of iterations
necessary in order to obtain
a dynamics which is far from that of an autonomous flow, has to grow faster than
$O(\delta^{-m})$, for every $m$. As we see, in order to obtain such kind of dynamics,
one has to control a very large number of iterations
of close to identity
maps, hence decompositions much longer than that given by (\ref{ruta}) have to be considered.

In this paper we propose such a decomposition (Theorem \ref{thm3}), using which we show that\\
{\em an arbitrary $C^r$-generic orientation-preserving $n$-dimensional
dynamics can be obtained by iterations of $C^r$-close to identity maps of $B^n$, $n\geq 2$.}

To make the formulations precise, we borrow some definitions
from \cite{T3}. Let $g$ be a $C^r$-diffeomorphism of a closed $n$-dimensional ball $D$.
Take an integer
$m>0$ and any $C^r$-diffeomorphism $\psi$ of $R^n$ such that $\psi(B^n)\subseteq D$.
The map $g_{m,\psi}=\psi^{-1}\circ g^m \circ \psi|_{B^n}$
is a $C^r$-diffeomorphism that maps $B^n$ into $R^n$. We will
call the maps $g_{m,\psi}$ obtained by this procedure {\em
renormalized iterations of $g$}.

\begin{theorem}$\!\!.$ \label{thm1} In the space of $C^r$-smooth orientation-preserving
diffeomorphisms of $B^n$ into $R^n$ ($n\geq 2$) there is a residual set ${\cal S}_r$ such that
for every map $F\in {\cal S}_r$, for every $\delta>0$ and for every $n$-dimensional ball $D$
there exists a map $g: R^n\rightarrow R^n$, equal to identity outside $D$,
such that $||g-id||_{C^r}<\delta$ and $F$ is a renormalized iteration of $g$.
\end{theorem}

In other words, a generic $C^r$-diffeomorphism $B^n\rightarrow R^n$ is,
up to a smooth coordinate transformation, an iteration of an arbitrarily close to identity map.
This theorem is proven in Section \ref{proof1}. Actually, we prove there that\\
{\em for any $\delta>0$ and $\varepsilon>0$, for
every orientation-preserving $C^r$-diffeomorphism $F:B^n\rightarrow R^n$
there exists a $\delta$-close to identity map $g$, equal to identity outside
a given ball $D$, such that $||F-g_{m,\psi}||_{C^r}<\varepsilon$ for
some renormalized iteration $g_{m,\psi}=\psi^{-1}\circ g^m \circ \psi|_{B^n}$}
(it is enough to prove this for one particular ball $D$ that may depend on $F$, but
the choice of $D$ should be fixed in advance, not depending on $\delta$ nor
on the accuracy $\varepsilon$ of the approximation; then for other balls the claim
will remain true because there always exists an affine conjugacy that takes one ball to the other).
Moreover, we construct the map $g$
and the coordinate transformation $\psi$ in such a way that the balls $g^i(\psi (B^n))$ ($i=0,\dots,m-1$)
do not
intersect each other --- hence, by adding small, localized in $g^{m-1}(\psi (B^n))$,
perturbations to $g$, we may make $g_{m,\psi}$ run an open subset in the space of $C^r$-diffeomorphisms
$B^n\rightarrow R^n$. Thus, the set ${\cal S}_r(\delta)$ of all renormalized iterations
of the maps $g: D\rightarrow D$
such that $||g-id||_{C^r}<\delta$ contains a subset which is open and dense in the space
of $C^r$-smooth orientation-preserving
diffeomorphisms of $B^n$ into $R^n$, for every $\delta>0$. Hence, the intersection ${\cal S}_r$
of these sets over all $\delta>0$ is residual (and independent
of the choice of $D$), which gives us the theorem.

The first step in our construction of the approximations of the given map $F$
by renormalized iterations
is a representation of $F$
as a superposition of a pair of certain special maps and some volume-preserving diffeomorphisms (Lemma
\ref{lm1}).
Each of the special maps can be realized as a flow map through a kind of saddle-node
bifurcation (see Fig.2), reminiscent of the so-called ``Iljashenko lips'' (see \cite{Ibook}).
For the volume-preserving diffeomorphisms one may adjust the results obtained in  \cite{T3}
for symplectic
diffeomorphisms and prove (Lemma \ref{lm2}) the existence of an arbitrarily good, in
the $C^r$-norm on any compact, polynomial approximation by a
superposition of volume-preserving H\'enon-like
maps. It is known that H\'enon-like maps often appear
as rescaled first-return maps near a homoclinic tangency
(cf. \cite{GTS3,GTS6,GTS7}). In this paper we find a kind of homoclinic tangency which
does incorporate all the H\'enon-like maps that appear in our volume-preserving
polynomial approximations.

Thus, we show that the map $F$ can be approximated
arbitrarily well by a superposition of maps related to certain
homoclinic bifurcations. The last step is to build a close to
identity map which displays these bifurcations simultaneously.
This is achieved by an arbitrarily small perturbation of
the time-$\delta$ map of a certain $C^\infty$ flow
(the time $\delta$ map of a flow is, obviously, $O(\delta)$-close to identity).

Note that the approximation (that we construct in Lemma \ref{lm2}) of any
volume-preserving diffeomorphism of a unit ball into $R^n$
by a polynomial volume-preserving diffeomorphism is not
straightforward, because the Jacobian of the approximating diffeomorphism
should be equal to $1$ everywhere, and this constrain is quite
strong for polynomial maps. Had the approximation result been
true for all volume-preserving maps, i.e. not necessarily
diffeomorphisms, it would produce a counterexample to the famous
``Jacobian conjecture''; however, our approximation uses in an essential way
the injectivity of the map that has to be approximated
(we represent the map as a shift by the orbits of some smooth non-autonomous flow).

It should be mentioned that Theorem \ref{thm1} does not hold true at $n=1$. Indeed,
if a map $F$ on the interval $B^1$ has two fixed points (with the multipliers different from $1$),
then every close map $\hat F$ has a pair of fixed points $P_{1,2}$ as well.
If such $\hat F$ is a renormalized iteration
of a diffeomorphism $g$, i.e. if $\hat F=\psi^{-1}\circ g^m \circ \psi$,
then $g$ will also have a pair of fixed points, $\psi (P_1)$ and $\psi (P_2)$
(at $n>1$ this is not true). The interval between $P_1$ and $P_2$ will therefore be invariant
with respect to $\psi^{-1}\circ g \circ \psi$, hence $\psi^{-1}\circ g \circ \psi$ will be
a root of degree $m>1$ of the map $\hat F$ on this interval. Now note that the maps of the interval
that have a root are not dense in $C^2$, according to \cite{Bel}. Thus we obtain that renormalized
iterations are not dense either.

One can check through the proof of Theorem \ref{thm1} that it holds
true for finite-parameter families of orientation-preserving diffeomorphisms:\\
{\em in the space of $k$-parameter
families $F_\varepsilon$, $\varepsilon\in B^k$, of $C^r$-smooth orientation-preserving
diffeomorphisms of $B^n$ into $R^n$ ($n\geq 2$) there is a residual set ${\cal S}_{kr}$ such that
for every $F_\varepsilon\in {\cal S}_{kr}$ and for every $\delta>0$ there exists
$g_\varepsilon: D\times B^k\rightarrow D$
such that $||g_\varepsilon-id||_{C^r}<\delta$ and
$F_\varepsilon=\psi_\varepsilon^{-1}\circ g_\varepsilon^m \circ \psi_\varepsilon|_{B^n}$
for some $m>0$ and some family $\psi_\varepsilon$ of $C^r$-diffeomorphisms of $R^n$.}\\

Thus, any dynamical phenomenon which occurs generically in finite-parameter families of dynamical
systems can be encountered in maps arbitrarily close to identity (with the same dimension of the phase space).

To put the result into a general perspective, we recall that one of
the main sources of difficulties in the theory of dynamical systems
is that structurally stable systems are not dense in the space of all systems
\cite{Sm,Ne74,Ne}, moreover most natural examples of chaotic dynamics are indeed
structurally unstable (like e.g. the famous Lorenz attractor \cite{ABS}). Understanding
the dynamics of systems from the open regions of structural instability
(i.e. the regions where arbitrarily close to every  system there
is a system which is not topologically conjugate to it) has
been the subject of active research for the past four decades. It often
happens, and helps a lot, that structurally unstable systems may
possess certain robust properties, i.e. dynamical properties
which are not destroyed by small perturbations. For example,
systems with Lorenz attractor are pseudohyperbolic (or volume-hyperbolic) \cite{TS98,TS7,Bon},
and this is, in fact, the very property which allowed for a very detailed description
of them \cite{Wi,ABS,Bu}. Another robust property is uniform partial hyperbolicity,
a rich theory of systems possessing it is actively developing \cite{Vi,BDV}.
In fact, not so much of robust properties are known, it could even happen that
beyond the mentioned partial hyperbolicity and volume-hyperbolicity
no other robust dynamical properties exist. This claim can be demonstrated for various
examples of homoclinic bifurcations (see \cite{T96}), and can be used as
a guiding principle in the study of bifurcations of systems with a non-trivial dynamics:\\
{\em given an $n$-dimensional system with a compact invariant set that
is neither partially- nor volume-hyperbolic,
every dynamics that is possible in $B^n$ should be expected to occur at the
bifurcations of this particular system.}\\
The last statement is not a theorem and it might
be not true in some situations, still it gives a useful view on global bifurcations.
In particular, it was explicitly applied in \cite{TZ} to Galerkin approximations of damped nonlinear wave
equations in order to obtain estimates from below on
the dimension of attractors in the situation
where classical methods \cite{BaV} do not work.

Theorem \ref{thm1} gives one more example to the above stated principle: the identity map has no kind of
hyperbolic structure, neither it contracts nor expands volumes, so it should not be surprising
that its bifurcations provide an ultimately rich dynamics.

The same idea can be expressed in somewhat different terms. Let us call
the set of all renormalized iterations of a map $g: D\rightarrow D$ its
{\em dynamical conjugacy class}. The map will be called {\em $C^r$-universal} (cf. \cite{T3})
if its dynamical conjugacy class is $C^r$-dense among all orientation-preserving
$C^r$-diffeomorphisms of the closed unit ball $B^n$ into $R^n$.
By the definition, the dynamics of any single universal map is
ultimately complicated and rich, and the detailed understanding of
it is not simpler than understanding of all diffeomorphisms $B^n\rightarrow R^n$ altogether.

At the first glance, the mere existence of $C^r$-univesal maps
of a closed ball is not obvious for
sufficiently large $r$. However, Theorem \ref{thm1} immediately implies the following

\begin{theorem}$\!\!.$ \label{thm2} For every $r\geq 1$, $C^r$-universal diffeomorphisms
of a given closed ball $D$ exist arbitrarily close, in the $C^r$-metric, to the identity map.
\end{theorem}

{\em Proof}. Take an arbitrary sequence of pairwise disjoint closed balls $D_j\subset D$,
a sequence of maps $F_j$ which is $C^r$-dense in space of orientation-preserving
$C^r$-diffeomorphisms $B^n\rightarrow R^n$, and a sequence $\varepsilon_j\rightarrow+0$
as $j\rightarrow +\infty$. By Theorem \ref{thm1}, given any $\delta$, there exist maps $g_j$ such that
$g_j$ is identity outside $D_j$, some renormalized iteration of $g_j$ is $\varepsilon_j$-close
to $F_j$, and $\|g_j-id\|_{C^r}\leq\delta$. By construction, the map $g(x)\equiv g_j(x)$ at $x\in D_j$
($j=1,2,\dots$) and $g(x)\equiv x$ at $x\in D\backslash \cup_{j=1}^{\infty} D_j$
is $C^r$-universal and $\delta$-close to identity.
~$\Box$

As an immediate consequence of Theorem \ref{thm2}, we note that\\
{\em arbitrarily close to identity map in the space of $C^r$-diffeomorphisms
$B^n\rightarrow B^n$ (any $r\geq1$, any $n\geq 2$) there exist maps with
infinitely many coexisting uniformly-hyperbolic attractors of all possible topological types.}\\
This is true because a hyperbolic attractor is a structurally stable object:
given a map with a uniformly-hyperbolic attractor, any $C^1$-close map has
a hyperbolic attractor topologically conjugate to the original one. For every $n\geq 2$ there
exists a $C^r$-diffeomorphism $B^n\rightarrow R^n$ with a hyperbolic attractor. Hence,
infinitely many of the maps $F_j$ in the proof of Theorem \ref{thm2} have a hyperbolic attractor
as well. It follows that each universal map constructed in Theorem \ref{thm2}
has infinitely many hyperbolic attractors (and repellers, for that matter).

Thus, we can further pursue the approach of \cite{RT} and claim that hyperbolic
attractors can be born at the third Andronov-Hopf bifurcation. Recall that at the
primary Andronov-Hopf bifurcation a limit cycle is born from an equilibrium state, and at the
secondary Andronov-Hopf bifurcation a two-dimensional invariant torus is born
from the limit cycle. The third Andronov-Hopf bifurcation occurs when a three-dimensional
invariant torus is born from the two-dimensional one (filled by quasiperiodic orbits).
The so-called Landau-Hopf scenario of the onset of turbulence envisioned a chain
of further Andronov-Hopf bifurcations which would lead to a creation of an
invariant torus of a sufficiently high dimension, i.e. to a quasiperiodic regime
with a high number of rationally independent frequencies (see more discussion in \cite{tbft}).
However, as Ruelle and Takens pointed out in \cite{RT}, the dynamics on
the invariant torus is not necessarily quasiperiodic: at the moment the torus is born,
the system can be perturbed in such a way that every orbit on the torus will become periodic,
hence the Poincar\'e map will be the identity, hence, as it follows from our results above,
a further small perturbation may lead to a chaotic dynamics for a flow
on the torus of dimension 3 or higher.

\section{Universal maps in the Newhouse domain}
\setcounter{equation}{0}
In general, it follows from Theorem \ref{thm2} that
every time we have a periodic orbit for which the corresponding
Poincar\'e map is, locally, identity:
$$\bar x =x,$$
or coincides with identity up to flat (i.e. sufficiently high order) terms:
$$\bar x= x + o(\|x\|^r),$$
a $C^r$-small perturbation of the system can make the Poincar\'e map universal,
i.e. bifurcations of this orbit can produce dynamics as complicated
as it only possible for the given dimension of the phase space. Thus, arbitrarily
complicated dynamical phenomena can be uncovered by the study of bifurcations
of periodic orbits alone.

In order to show how powerful this observation can be, let us apply it to the analysis
of the dynamics of two-dimensional diffeomorphisms from the Newhouse domain. In the space of
$C^r$-smooth dynamical systems on any smooth manifold, Newhouse
domain is the interior of the closure of the set of systems that have a homoclinic tangency
(a tangency between the stable and unstable manifolds of a saddle periodic orbit; for a saddle periodic
orbit to exist, the dimension of the phase space should be at least $3$ in the case of continuous time
and at least $2$ for discrete dynamical systems). A non-trivial fact
\cite{Ne74,Ne,GST93,PV94}
is that the Newhouse domain at $r\geq 2$ is always non-empty
and adjoins to every system with a homoclinic tangency. Importantly, most of known global bifurcations
which lead to the emergence of chaotic dynamics or happen within the class of systems
with complex (chaotic) behavior are accompanied by a creation of homoclinic tangencies. Therefore,
Newhouse regions in the space of parameters can be detected for virtually every dynamical model with
chaos (see more discussion in \cite{GST92,Y,T3,GTS6,GTS7}). In my opinion, the ubiquitous presence of
homoclinic tangencies in the dynamical models of a natural origin makes the study of maps
from the Newhouse domain the most important problem of chaotic dynamics. It should also be mentioned that
a commonly shared believe (with no hope to prove yet) is that the space of two-dimensional
$C^r$-diffeomorphisms with $r\geq 2$ is the closure of the union of just three open sets: Morse-Smale
systems (i.e. systems with simple dynamics), axiom A systems, and the Newhouse domain, i.e. unless
a two-dimensional map with a chaotic dynamics is uniformly hyperbolic, it most probably lies in the
Newhouse domain.

Typically, a map from the Newhouse domain possesses an invariant basic hyperbolic set which
is wild -- the term meaning that for the map itself, and for every $C^r$-close map ($r\geq 2$),
there exists a pair of orbits within the hyperbolic set such that the unstable manifold of one orbit
has a quadratic tangency with the stable manifold of the other \cite{Ne74,Ne}. For a fixed pair of orbits,
the corresponding tangency
is a codimension-$1$ bifurcation phenomenon, so it can always be removed by a small smooth perturbation;
the wildness, nevertheless, means that once the original tangency is removed a new tangency
always appears,
corresponding to some other pair of orbits from the same hyperbolic set.

A wild hyperbolic set of a $C^r$-diffeomorphism of a two-dimensional smooth manifold will be called
{\em ultimately wild} if it contains a pair of periodic orbits such that the saddle value
at one periodic orbit is less than $1$ and at the other periodic orbit
it is greater than $1$. The open subset of the Newhouse domain which corresponds to maps with ultimately
wild hyperbolic sets will be called the {\em absolute Newhouse domain}. The saddle value is, by definition,
the absolute value of the product of multipliers of the periodic orbit. If $Q$ is a period $q$ point of
a map $f$ (i.e. $f^q Q=Q$), its multipliers are the eigenvalues of the matrix of derivatives of $f^q$
at the point $Q$. Therefore, if the saddle value is greater than $1$ in the absolute value,
then the map $f$ expands area near $Q$, and it is area-contracting near $Q$ if the saddle value is less
than $1$. So, no map from the absolute Newhouse domain is area-contracting, nor area-expanding. Moreover,
the persistent tangencies between the stable and unstable manifolds of the wild hyperbolic set mean that
none of these maps is uniformly hyperbolic, nor uniformly partially-hyperbolic. Thus, there is no
obvious restrictions on the dynamics of two-dimensional diffeomorphisms from the absolute Newhouse
domain, and the following theorem is, therefore, in agreement with the general ``guiding principle''
formulated in Section \ref{rtu}.

\begin{theorem}$\!\!.$ \label{thmn} For every $r\geq 2$, a $C^r$-generic
diffeomorphism from the absolute Newhouse domain is $C^r$-universal.
\end{theorem}

{\em Proof}. Fix any $r\geq 2$ and let a $C^r$-diffeomorphism of a smooth two-dimensional manifold possess
an ultimately wild basic hyperbolic
set $\Lambda$. Let $P$ and $Q$ be two saddle periodic points in $\Lambda$ ($f^pP=P$ and $f^qQ=Q$) with the
saddle values $\sigma_{_P}$ and $\sigma_{_Q}$ such that $\sigma_{_P}<1$ and $\sigma_{_Q}>1$.
As the periodic points $P$ and $Q$ are hyperbolic, they are preserved at all sufficiently
small perturbations. Recall also that
the unstable manifold of every point in the basic hyperbolic set has a transverse
intersection with the stable manifold of every other orbit in this set. Therefore, the invariant unstable
manifold $W^u(P)$ intersects the invariant stable manifold $W^s(Q)$ transversely
at the points of some heteroclinic orbit
$\Gamma_{PQ}$. By the transversality, this orbit is preserved for all maps sufficiently close to $f$.
According to \cite{GTS6}, given every $m\geq 1$, in any neighborhood of $f$ in the $C^r$-topology,
there is a $C^\infty$-diffeomoprphism for which the unstable manifold $W^u(Q)$ has a tangency of order
$m$ to the stable manifold $W^s(P)$ at the points of some heteroclinic orbit $\Gamma_{QP}$ (Fig.3). This
is a non-trivial statement: while the possibility to obtain a quadratic tangency (i.e. $m=1$) by
an arbitrarily small perturbation follows immediately from the wildness of $\Lambda$
and from the fact that the stable and unstable manifold of any given periodic orbit in the basic
hyperbolic set are dense within the union of stable and, respectively, unstable
manifolds of all orbits in this set \cite{Ne74,Ne}, the tangencies of higher orders are due to
the existence of moduli of local $\Omega$-conjugacy \cite{GST92,GTS93b,GTS1}.

Let $\tilde f\in C^\infty$ be a $C^r$-close to $f$ diffeomorphism for which the above described
heteroclinic cycle exists.
This cycle is a closed set that consists of $4$ orbits: two periodic orbits (the orbits of $P$
and $Q$) with $\sigma_{_P}<1$ and $\sigma_{_Q}>1$, a transverse heteroclinic $\Gamma_{PQ}$
and the orbit $\Gamma_{QP}$ of heteroclinic tangency of
a sufficiently high order $m$. We will call an open region filled by periodic orbits of the
same period {\em a periodic spot}. In Section \ref{henhom} (see Lemma \ref{thmn2} and Remark after it)
we show that\\
{\em for any neighborhood $U$ of the heteroclinic cycle, arbitrarily $C^r$-close to $\tilde f$ there is
a diffeomorphism $\hat f$ which has a periodic spot whose all iterations lie in $U$.}\\
For every orbit in the periodic spot the corresponding
Poincar\'e map (the map $\hat f^k$ where $k$ is the period of the spot)
is, locally, identity $\bar x =x$. Hence, as we mentioned in the beginning of this Section,
it follows from Theorem \ref{thm2} that in any $C^r$-neighborhood
of $\hat f$ there exist $C^r$-universal maps.

We have just shown (modulo Lemma \ref{thmn2} which is proved in Section \ref{henhom})
that universal maps are dense in the absolute Newhouse domain. This immediately
implies the genericity of the universal maps. Indeed, given an orientation-preserving $C^r$-diffeomorphism
$g$ of the unit disk into $R^2$, denote as ${\cal V}(g,\delta)$ the set of
all $C^r$-diffeomorphisms from the absolute Newhouse domain whose dynamical
conjugacy class intersects the open $\delta$-neighborhood of $g$
(i.e. whose certain renormalized iteration is at a $C^r$-distance smaller than
$\delta$ from $g$). This set, by definition, contains
all universal maps --- hence, it is dense in the absolute Newhouse domain. This set is also
open by definition.
Take a countable sequence of maps $g_i$ which is dense
in the space of orientation-preserving $C^r$-diffeomorphisms of the unit disk into $R^2$, and a
sequence $\delta_j$ of positive reals converging to zero. By construction,
the countable intersection $\cap_{i,j} {\cal V}(g_i,\delta_j)$ is a residual subset
of the absolute Newhouse domain, and every map that belongs to this set is universal.
~$\Box$

In essence, this theorem gives an exhaustive characterization of the richness of dynamical
behavior in the absolute Newhouse domain: every two-dimensional dynamics can
be approximated by iterations of any generic map from this domain. In a simpler case
of the Newhouse domain in the space of area-preserving maps,
a similar statement is contained in \cite{GTS7}: iterations of a generic
area-preserving map from the Newhouse domain approximate all symplectic dynamics in a two-dimensional
disc. For {\em area-contracting} maps, it follows from \cite{GTS1} that the closure of the
dynamical conjugacy class of a generic map from the Newhouse domain contains all
{\em one-dimensional} maps (we cannot have truly two-dimensional dynamics here, as the areas are
contracted). As we mentioned, any two-dimensional map which is not in the class of maps
with a uniformly-hyperbolic
structure, nor on the boundary of this class, falls, hypothetically,
in one of the three types of the Newhouse domain: the first is
filled by area-contracting maps, the second by area-expanding maps (i.e. the maps inverse to the maps
of the first type), and the third is the absolute Newhouse domain. We see that our Theorem \ref{thmn}
somehow completes the description of two-dimensional dynamics.

The three types of Newhouse domains of two-dimensional maps were introduced in \cite{GTS97}. It has been
known since \cite{Ne74} that a generic area-contracting map from the Newhouse domain
has an infinite set of stable periodic orbits, and the closure of this set contains a (wild) hyperbolic
set. The latter fact is especially important: while chaotic dynamics is usually associated with hyperbolic
sets, i.e. with saddle orbits, the Newhouse result shows that stable periodic motions can imitate chaos
arbitrarily well, and they indeed do it generically. In \cite{GTS97}, for the Newhouse domain of the
third type, it was shown that a generic map has
both an infinite set of stable periodic orbits and an infinite set of repelling periodic orbits,
moreover the intersection of these sets is non-empty and contains an ultimately-wild hyperbolic set.
Thus, not only the Newhouse phenomenon holds, we also have a new effect here: a generic inseparability
of attractors from repellers. In fact, these attractors and repellers can
be more complicated than just periodic orbits: in \cite{St}, the coexistence of infinitely many
closed invariant curves was established for two-dimensional maps from
the Newhouse domains of the third type. Our Theorem \ref{thmn} strengthens these observations:
it implies the coexistence of infinitely many hyperbolic attractors and repellers for a generic
map from the absolute Newhouse domain (see remarks to Theorem \ref{thm2}). As we obtain the hyperbolic
attractors and repellers by a perturbation of periodic spots, and the periodic spots are found in
an arbitrarily small neighborhood of any heteroclinic cycle of the type we consider in
Theorem \ref{thmn}, it follows that the closures of the set of hyperbolic attractors and of the set
of hyperbolic repellers that we construct here contain any transverse heteroclinic orbit connecting
the points $P$ and $Q$. Such heteroclinic orbits are dense in the basic hyperbolic set $\Lambda$. Thus,
it follows from the proof of the theorem that\\
{\em a $C^r$-generic two-dimensional map from the absolute Newhouse domain has an infinite set
of uniformly-hyperbolic attractors and an infinite set of uniformly-hyperbolic repellers, and
the intersection of the closures of these sets contains a non-trivial hyperbolic set.}\\

In what follows we prove Theorem 1 and finish the proof of Theorem \ref{thmn}.

\section{An approximation theorem.}
\setcounter{equation}{0}
Let $F$ be an orientation-preserving $C^r$-diffeomorphism ($r\geq 3$)
which maps the ball $B^n:\{\sum_{i=1}^n x_i^2\leq 1\}$ into $R^n$.
Without loss of generality we may assume that $F$ is extended onto the whole $R^n$, i.e. it becomes
a $C^r$-diffeomorphism $R^n\rightarrow R^n$, and it is identical (i.e. $F(x)=x$) at $\|x\|$ sufficiently
large; such extension is always possible. Let $K$ be a constant such that
\begin{equation}\label{jbound}
\sup_{x\in R^n} \frac{\|\nabla J(x)\|}{J(x)} < K,
\end{equation}
where $J(x)$ is the Jacobian of $F$. Denote $R^n_+:=\{x_n>0\}$.

\begin{lemma}$\!\!.$ \label{lm1} There exists a pair of volume-preserving, orientation-preserving
$C^{r-2}$-diffeomorphisms $\Phi_1:R^n_+\rightarrow R^n_+$ and $\Phi_2:R^n\rightarrow R^n$
such that
\begin{equation}\label{decom}
F=\Phi_2\circ\Psi_2\circ\Phi_1\circ\Psi_1,
\end{equation}
where
\begin{equation}\label{psi12}
\Psi_j:=(x_1,\dots,x_{n-1},x_n)\mapsto (x_1,\dots,x_{n-1},\psi_j(x_n)) \;\;\;\;(j=1,2),
\end{equation}
with
\begin{equation}\label{phijs}
\psi_1(x_n)=e^{Kx_n}, \qquad \psi_2(x_n)=\ln x_n.
\end{equation}
\end{lemma}

{\em Proof}. We need to construct a volume-preserving diffeomorphism
$\Phi_1:(x_1,\dots,x_n\geq 0)\mapsto (\bar x_1,\dots,\bar x_n\geq 0)$ in such a way that
\begin{equation}\label{req}
\det\frac{\partial}{\partial x}\Psi_2\circ\Phi_1\circ\Psi_1(x)\equiv J(x)
\end{equation}
(then the Jacobian of $\Phi_2=F\circ(\Psi_2\circ\Phi_1\circ\Psi_1)^{-1}$ will be equal to $1$
automatically). By (\ref{psi12}),(\ref{phijs}), condition (\ref{req}) is equivalent to
$$\bar x_n=\phi(x_1,\dots,x_n)\equiv
\frac{K x_n}{J(x_1,\dots,x_{n-1},\frac{1}{K}\ln x_n)}.$$

It follows from (\ref{jbound}) that
\begin{equation}\label{evrw}
\partial \phi/\partial x_n>0
\end{equation}
everywhere. Moreover,
as $F$ is the identity map outside a bounded region of $R^n$, we have that
\begin{equation}\label{phikxn}
\phi(x)=K x_n
\end{equation}
outside a compact subregion of $R^n_+$. Therefore, every trajectory of the vector field
\begin{equation}\label{vf}
\dot x_j=0 \;\; (j\leq n-2), \;\;\; \dot x_{n-1}=\frac{\partial \phi}{\partial x_n}, \;\;\;
\dot x_{n}=-\frac{\partial \phi}{\partial x_{n-1}}
\end{equation}
extends for all $x_{n-1}\in(-\infty,+\infty)$, and the time $\tau(x)$ that the trajectory
of the point $x$ needs to get to $x_{n-1}=0$ is a $C^{r-2}$-function of $x$, well defined everywhere in $R^n_+$. Moreover, as it follows from (\ref{phikxn}),(\ref{vf})
\begin{equation}\label{taukxn}
\tau(x)=- \frac{1}{K} x_{n-1}+\tau_0(x),
\end{equation}
where $\tau_0$ is a uniformly bounded function, vanishing identically at $x_n$ close
to zero and at sufficiently large $x_n$. Thus, for every fixed values of $x_j$
($j\leq n-2$), given any $C\in(-\infty,+\infty)$ the level line $\tau(x)=C$ in the $(x_{n-1},x_n)$-plane coincides with the straight line $x_{n-1}=-KC$ at $x_n$ close
to zero and at sufficiently large $x_n$. Every such level line is a connected set
(as it is the image of the line $\{x_{n-1}=0, x_n\geq 0\}$ by the
time-$(-C)$ map of the flow of (\ref{vf})). Thus, as $x_n$ runs from $0$ to $+\infty$,
the value of $\phi$ on this line runs all the values from $0$ to $+\infty$ (see (\ref{phikxn})). It follows that the map $R^n_+\rightarrow R^n_+$ defined by
\begin{equation}\label{phi2}
\bar x_j=x_j \;\; (j\leq n-2), \;\;\; \bar x_{n-1}=-\tau(x), \;\;\; \bar x_{n}=\phi(x)
\end{equation}
is surjective.

By (\ref{vf}),
\begin{equation}\label{dettf}
\frac{\partial \tau}{\partial x_{n-1}}\; \frac{\partial \phi}{\partial x_n}-
\frac{\partial \tau}{\partial x_n}\; \frac{\partial \phi}{\partial x_{n-1}}=-1.
\end{equation}
It follows that for every fixed values of $x_j$ ($j\leq n-2$), the function $\phi$ changes monotonically along every level line of $\tau$, which implies the injectivity of
map (\ref{phi2}). Thus, map (\ref{phi2}) is a $C^{r-2}$-diffeomorphism
$R^n_+\rightarrow R^n_+$. By (\ref{dettf}), it is volume-preserving and orientation-preserving, i.e. it is the sought map $\Phi_1$.
~$\Box$

The maps $x\mapsto\bar x$ of the following form:
\begin{equation}\label{hvol}
\bar x_1=x_2,\; \dots,\; \bar x_{n-1}=x_n, \quad \bar x_n = (-1)^{n+1} x_1 + h(x_2,\dots,x_n)
\end{equation}
(note no dependence on $x_1$ in $h$), will be called {\em H\'enon-like volume-preserving maps}. Note that such maps are always one-to-one,
and the inverse map is also H\'enon-like.

\begin{theorem}$\!\!.$ \label{thm3} Every orientation-preserving $C^r$-diffeomorphism
$F: B^n\rightarrow R^n$ can be arbitrarily closely approximated, in the $C^r$-norm on $B^n$,
by a map of the following form:
\begin{equation}\label{hdec}
H_{2q_2}\circ\dots\circ H_{21}\circ\Psi_2\circ H_{1q_1}\circ \dots \circ H_{11}\circ \Psi_1,
\end{equation}
where the maps $\Psi_{1,2}$ are given by (\ref{psi12}), and $H_{js}$ ($j=1,2; s=1,\dots q_j$) are certain
polynomial H\'enon-like volume-preserving maps.
\end{theorem}

{\em Proof}. First, take a $C^{r+2}$-diffeomorphism $\hat F$ which approximates $F$ sufficiently
closely in $C^r$. For the map $\hat F$ construct the decomposition
$\hat F=\Phi_2\circ\Psi_2\circ\Phi_1\circ\Psi_1$ given by Lemma \ref{lm1}; all
the maps in the decomposition are at least $C^r$. The map $\Phi_1$
can be extended onto $x_n\leq 0$ by the rule
$\bar x_n= Kx_n, \;\bar x_{n-1}=x_{n-1}/K$ (see (\ref{phi2}),(\ref{phikxn}),(\ref{taukxn})),
so it becomes a volume-preserving, orientation-preserving $C^r$-diffeomorphism
$R^n\rightarrow R^n$. Then the theorem follows immediately from Lemma \ref{lm2} below.
~$\Box$

\begin{lemma}$\!\!.$ \label{lm2} Every volume-preserving, orientation-preserving $C^r$-diffeomorphism
$\Phi:R^n\rightarrow R^n$ can be arbitrarily closely approximated, in the $C^r$-norm on any
given compact, by a composition of polynomial H\'enon-like volume-preserving maps.
\end{lemma}

{\em Proof}. At $n=2$ this result is immediately given by Theorem 1 in \cite{T3},
so we proceed to the case $n\geq 3$. It is well known that $\Phi$ can be imbedded in a smooth
in $t$ family ${\cal F}_{t}$
of volume-preserving $C^r$-diffeomorphisms $R^n\rightarrow R^n$ such that ${\cal F}_{0}\equiv id$
and ${\cal F}_{1}=\Phi$. The derivative $\frac{d}{dt}{\cal F}_{t}$ defines a divergence-free
vector field $X(t,x)$, i.e. the diffeomorphism ${\cal F}_{t}$ is the time-$t$ shift
by the flow generated by the field $X$. One can approximate $X$ arbitrarily closely on any given
compact by a $C^\infty$-smooth divergence-free vector field which is defined and bounded
for all $(x,t)\in R^n\times [0,1]$.
Therefore, it is enough to prove the lemma only for those $\Phi$ which can be obtained as the time-$1$
shift by the flow generated by such a vector field, i.e. we may assume that $X\in C^{\infty}_b$
with no loss of generality.

Let $T_{\tau t}={\cal F}_{t+\tau}\circ {\cal F}_{t}^{-1}$, i.e. it is the shift by the flow of
$X$ from the time $t$ to $t+\tau$. This map is $O(\tau)$-close to identity, in the $C^r$-norm
on any compact subset of $R^n\times [0,1]$. By construction, given any arbitrarily large integer $N$,
\begin{equation}\label{fdeco}
\Phi=T_{\tau, (N-1)\tau}\circ\dots\circ  T_{\tau, m\tau}\circ\dots\circ T_{\tau, 0}
\end{equation}
where $\tau=1/N$, and $m=0,\dots, N-1$.

Note that the vector field $X$ admits the following representation:
\begin{equation}\label{sumfi} X=\sum_{i=1}^{n-1} X^{(i)}\end{equation}
where $X^{(i)}$ is a $C^\infty$-smooth divergence-free vector field such that
\begin{equation}\label{hamf}
\dot x_j\equiv 0 \;\;\mbox{ at }\;\; j\neq i, i+1.
\end{equation}
Indeed, if we write the field $X$ as
$$\dot x_i = \xi_i(x,t), \quad i=1,\dots, n,$$
where $\displaystyle \sum_{i=1}^n \frac{\partial \xi_i}{\partial x_i}\equiv 0$ (the
zero divergence condition),
then the fields $X^{(i)}$ are defined as
$$\dot x_i=\eta_i(x,t), \;\;\; \dot x_{i+1}=\zeta_i(x,t)$$
with
$$\eta_1\equiv \xi_1, \quad \eta_i\equiv\xi_i-\zeta_{i-1}\;\; (i=2,\dots,n-1),$$
$$\zeta_i=-\int_0^{x_{i+1}}
\frac{\partial}{\partial x_i}\eta_i(x_1,\dots,x_i,s,x_{i+2},\dots,x_n,t) ds \;\;
(i=1,\dots,n-2), \quad \zeta_{n-1}\equiv\xi_n.$$
By construction, the fields $X^{(1)},\dots,X^{(n-2)}$ are divergence-free, and
$X^{(n-1)}=X - X^{(1)}-\dots-X^{(n-2)}$, so $X^{(n-1)}$ is also divergence-free, as
$X$ is.

It follows from (\ref{sumfi}) that
\begin{equation}\label{prodf}
T_{\tau t}=T_{\tau t}^{(n-1)}\circ\dots\circ T_{\tau t}^{(1)}+O(\tau^2),
\end{equation}
where $T_{\tau t}^{(i)}$ is the shift by the flow generated by the vector field $X^{(i)}$.
Recall that the maps $T_{\tau, i\tau}$ in (\ref{fdeco}) are $O(1/N)$-close to identity.
Therefore, it follows from (\ref{prodf}),(\ref{fdeco}) that
\begin{equation}\label{fdecoa}
\Phi=T_{\tau,(N-1)\tau}^{(n-1)}\circ\dots\circ T_{\tau,(N-1)\tau}^{(1)}\circ\dots\circ
T_{\tau,m\tau}^{(n-1)}\circ\dots\circ T_{\tau,m\tau}^{(1)}\circ\dots\circ
T_{\tau,0}^{(n-1)}\circ\dots\circ T_{\tau,0}^{(1)}+O(\tau),
\end{equation}
uniformly on compacta.

As $\tau$ can be taken arbitrarily small, it follows that in order to prove
the lemma, it is enough to prove that every of the maps $T^{(i)}_{\tau t}$
in the right-hand side of
(\ref{fdecoa}) can be approximated arbitrarily well by a composition of H\'enon-like volume-preserving maps.
The maps $T^{(i)}_{\tau t}$ are volume-preserving and satisfy
\begin{equation}\label{fixx}\bar x_j=x_j \;\;\mbox{ at }\;\; j\neq i, i+1\end{equation}
(see (\ref{hamf})). Therefore, if we denote
\begin{equation}\label{titat} \bar x_i=p(x),\quad \bar x_{i+1}=q(x),\end{equation}
then
\begin{equation}\label{2dyac}\frac{\partial(p,q)}{\partial(x_i,x_{i+1})}=1.\end{equation}
Thus, we can view (\ref{titat}) as an $(n-2)$-parameter family of symplectic two-dimensional maps
$(x_i,x_{i+1})\mapsto (\bar x_i,\bar x_{i+1})$ parametrized by
$(x_1,\dots,x_{i-1},x_{i+2},\dots,x_n)$.

According to \cite{T3}, every finite-parameter family of symplectic maps can be approximated (on any
compact) by a composition of families of H\'enon-like maps, i.e., in our case, maps of the form
$$\bar x_i=x_{i+1}, \quad \bar x_{i+1}=-x_i + h(x_{i+1};x_1,\dots,x_{i-1},x_{i+2},\dots,x_n).$$
It follows that every map of the form (\ref{fixx}),(\ref{titat}),(\ref{2dyac}) can
be approximated arbitrarily closely by a composition of the maps of the form
\begin{equation}\label{semihen}
\begin{array}{l}
\bar x_j=x_j \;\;\mbox{ at }\;\; j\neq i, i+1,\\
\bar x_i=x_{i+1},\\
\bar x_{i+1}=-x_i + h(x_{i+1};x_1,\dots,x_{i-1};x_{i+2},\dots,x_n).
\end{array}
\end{equation}
It just remains to note that
every map of form (\ref{semihen}) is a composition of volume-preserving
H\'enon-like maps; namely, it equals to
$$S^{n-i-1}\circ H\circ S\circ Q^{n-1}\circ S^{i+1},$$
where
$$S:=(x_1,\dots,x_n)\mapsto (x_2,\dots,x_n,(-1)^{n+1}x_1),$$
$$Q:=(x_1,\dots,x_n)\mapsto (x_2,\dots,x_n,\sum_{j=1}^n (-1)^{n+j}x_j),$$
$$\begin{array}{l}
H:=\{\bar x_1=x_2,\dots,\bar x_{n-1}=x_n, \; \bar x_n = \\ \displaystyle \quad
= \;\sum_{j=1}^{n-1} (-1)^{n+j}x_j - x_n +h(x_n;x_{n-i+1},\dots,x_{n-1};
(-1)^{n+1}x_2,\dots,(-1)^{n+1}x_{n-i})\}.
\end{array}$$
End of the proof. ~$\Box$

\noindent Remark. Consider the map
\begin{equation}\label{phi0}
\Phi_0:=(x_1,\dots,x_{n-2},x_{n-1},x_n)\mapsto(x_1,\dots,x_{n-2},x_n,-x_{n-1}).
\end{equation}
This is an orientation- and volume- preserving diffeomorphism of $R^n$, therefore we may rewrite
(\ref{decom}) as follows:
$$F=\Phi_0\circ\tilde \Phi_2\circ\Psi_2\circ\Phi_0\circ\tilde\Phi_1\circ\Psi_1,$$
where $\tilde\Phi_{1,2}$ are orientation-preserving, volume-preserving
$C^{r-2}$-diffeomorphisms
($\tilde\Phi_j=\Phi_0^{-1}\circ\Phi_j$; we assume that $\Phi_1$ is extended onto the whole of $R^n$,
like in Theorem \ref{thm3}). Now, by Lemma \ref{lm2}, we obtain the following, more
convenient for us, analog of Theorem \ref{thm3}:\\ {\em the map $F$
can be arbitrarily closely approximated by a map of the following form:
\begin{equation}\label{dec0n}
\Phi_0\circ \tilde H_{2q_2}\circ \dots \circ\tilde H_{21}\circ\Psi_2\circ \Phi_0\circ
\tilde H_{1q_1}\circ \dots \circ \tilde H_{11}\circ \Psi_1,
\end{equation}
with polynomial H\'enon-like volume-preserving maps $\tilde H_{js}$.}

\section{Proof of Theorem 1.}\label{proof1}

\renewcommand{\theequation}{\thesubsection.\arabic{equation}}

Given a diffeomorphism $F:B^n\rightarrow R^n$ we will take its sufficiently
close approximation $\hat F$ in the form of (\ref{dec0n}).
Then we will construct a close to identity map (which we denote $\tilde Y_{\delta}$ below)
whose some renormalized iteration is a close (as close as we want)
approximation to $\hat F$.

The map $\tilde Y_{\delta}$ is a small perturbation (as small as we want, in
the $C^r$-norm for any aforehand given $r$) of the time-$\delta$ map $Y_\delta$
of a certain $C^\infty$ flow $Y$ in $R^n$; the constant $\delta$ can be chosen as small as we need
(i.e. $Y_\delta$ is indeed a small perturbation of the identity). In our construction,
the vector field of $Y$ vanishes identically
outside some ball $D$ that does not depend on the choice of the approximation $\hat F$,
it does not depend on $\delta$ either.
The small perturbations which we will apply to $Y_\delta$ will also be localized
in $D$. Hence, our close to identity maps $\tilde Y_\delta$ are all equal to identity
outside the same ball $D$. Since their renormilzed iterations approximate $F$ arbitrarily
well, this gives us Theorem 1, as we explained it in the remark after the theorem.

We define the flow $Y$ by means of the following procedure: we give explicit formulas
for the vector field inside certain blocks $U_{1\pm}, U_{2\pm}, V_{1,2}$ described below,
while between the blocks we specify only the transition time from the boundary
of one block to another and the corresponding Poincar\'e map.
The existence of a $C^\infty$ flow with arbitrary (of class $C^\infty$)
transition times and orientation-preserving Poincar\'e maps between block boundaries
is a routine fact (at least for the given geometry of the blocks, see Fig.2).

The idea of the construction is as follows.
Given a diffeomorphism $F:B^n\rightarrow R^n$, we may always approximate it by a
$C^\infty$-diffeomorphism arbitrarily well, so we assume that $F\in C^\infty$ from the
beginning. Let $\Phi_{1,2}$ and $\Psi_{1,2}$
be the maps defined by decomposition (\ref{decom}) of $F$. We define the vector field inside
the blocks $U_{1\pm}, U_{2\pm}$ in such a way that a kind of saddle-node bifurcation
is created inside each of the blocks (see (\ref{yw2+}),(\ref{yw+}); we build very
degenerate saddle-nodes in order to make formulas for the time-$t$ map simpler --
see (\ref{tmap2}),(\ref{tmap})). The Poincar\'e maps from the boundary
$\Sigma_{j+}^{out}$ of $U_{j+}$ to the boundary $\Sigma_{j-}^{in}$ of $U_{j-}$ ($j=1,2$)
are chosen in such a way that the resulting flow map from entering $U_{j+}$ till exiting $U_{j-}$
is linearly conjugate to the map $\Psi_j$
(see (\ref{stmap2}),(\ref{stmap}),(\ref{psi12}),(\ref{phijs})).

We make the flow in $V_{1,2}$ volume-preserving and linear (see (\ref{yw02}),
(\ref{yw0}),(\ref{gam})). Moreover, we put saddle equilibria into $V_{1,2}$.
We define the Poincar\'e map
between the boundaries $\Pi_{j+}^{out}$ and $\Pi_{j-}^{in}$ of $V_j$ (see Fig.2)
in such a way that a homoclinic loop to the saddle is created. For the time-$\delta$ map
$Y_\delta$ of the flow $Y$ the saddle equilibrium is a saddle fixed point,
and the homoclinic loop is a continuous family of homoclinic orbits. We perturb
the map $Y_\delta$ in such a way that this family splits into a finite set
of orbits of homoclinic tangency of sufficiently high orders and unfold these tangencies
then. The exact form of the perturbation (see (\ref{pcr2}),(\ref{pcr})) may be chosen such
that the iteration of the perturbed map $\tilde Y_\delta$ that corresponds to one round
near the homoclinic loop is a close approximation (in appropriately scaled coordinates)
to any given polynomial conservative H\'enon-like map (see (\ref{scsdl2}),(\ref{scsdl})).
Hence, a multi-round iteration of $\tilde Y_\delta$ can be made arbitrarily close to
a superposition of a finite number of any given polynomial H\'enon-like maps. By
Lemma \ref{lm2}
such superpositions approximate any given volume-preserving maps, e.g. maps $\Phi_j$
from the decomposition (\ref{decom}). In this way an iteration of the map
$\tilde Y_\delta$ that, after a large number of rounds near homoclinic orbits,
takes points entering $V_j$ to the points entering $U_{3-j}$ is made as
close as we want to the map $\Phi_j$ ($j=1,2$), in some rescaled coordinates.

The parameters of the flow and perturbations are chosen in such a way that
the coordinate scalings we make in the blocks $U_{1,2}$ and $V_{1,2}$ match each other.
Thus, by construction, the (renormalized) iteration of $\tilde Y_\delta$ that corresponds
to passing from the entrance to $U_{1+}$ through $U_{1-}$ into $V_1$, then to many rounds
near the homoclinic loop, then to exiting $V_1$ and entering $U_{2+}$, passing through
$U_{2-}$ to $V_2$, a number of near-homoclinic rounds and return to $U_{1+}$, is a close
approximation to $\Phi_2\circ\Psi_2\circ\Phi_1\circ\Psi_1$, i.e. to the original
diffeomorphism $F$ indeed.

In the two-dimensional case, the key fact that every area-preserving
diffeomorphism can be approximated by some multi-round iteration of a map with
a homoclinic tangency was proven in \cite{GTS6}. To deal with dimensions higher
than 2, we construct a very degenerate homoclinic tangency (in terms of \cite{T96},
both critical and Lyapunov dimensions for this tangency are equal to $n$ -- this
is necessary in order not to lose dimension at the rescaling of the first-return map).
We do not undertake an analysis of the corresponding bifurcation; instead, we
make explicit computations of the rescaled return maps for our particular example only.

\subsection{Two-dimensional case}
\setcounter{equation}{0}
To make computations more transparent we start with the case $n=2$.
Let $\Phi_{1,2}$ and $\Psi_{1,2}$ be the maps defined by (\ref{decom}). Let $I_{1\pm}$ and
$I_{2\pm}$ be intervals of values of $x_2$ such that $x_2\in I_{1+}$ at $(x_1,x_2)\in B^2$,
$x_2\in I_{1-}$ at $(x_1,x_2)\in \Psi_1(B^2)$, $x_2\in I_{2+}$ at
$(x_1,x_2)\in \Phi_1\circ\Psi_1(B^2)$
and $x_2\in I_{2-}$ at $(x_1,x_2)\in \Psi_2\circ\Phi_1\circ\Psi_1(B^2)$.
Let $R$ be such that all the intervals $I_{j\pm}$ lie within $\{|x_2|\leq R\}$.
Choose numbers $a_{1+}=a_{1-}+3=b_1+6=a_{2+}+9=a_{2-}+12=b_2+15$. Let the vector field of
$Y$ in the regions
$U_{j\sigma}:\{|x_1-a_{j\sigma}\|\leq 1, |x_2|\leq R\}$,
$j=1,2$, $\sigma=\pm 1$, be equal to
\begin{equation}\label{yw2+}
\begin{array}{l}\displaystyle
\dot x_1=-\mu_j-(1-\mu_j)(1-\xi(x_1-a_{j\sigma})),\\ \displaystyle
\dot x_2=\sigma \gamma_j \; x_2 \;\xi(x_1-a_{j\sigma}),
\end{array}
\end{equation}
where $\mu_{1,2}>0$ are small (see (\ref{mdf2})), $\gamma_{1,2}\in[0,1]$
(see (\ref{mdfa2})), and $\xi$ is a $C^\infty$ function such that
\begin{equation}\label{lam}
0\leq\xi\leq 1, \;\;\; \xi(0)=1, \;\;\; \xi(z)\equiv 0 \; \mbox{ at } \;
|z|\geq \frac{1}{2}.
\end{equation}

\begin{figure}
\centerline{\includegraphics[height=15.5cm]{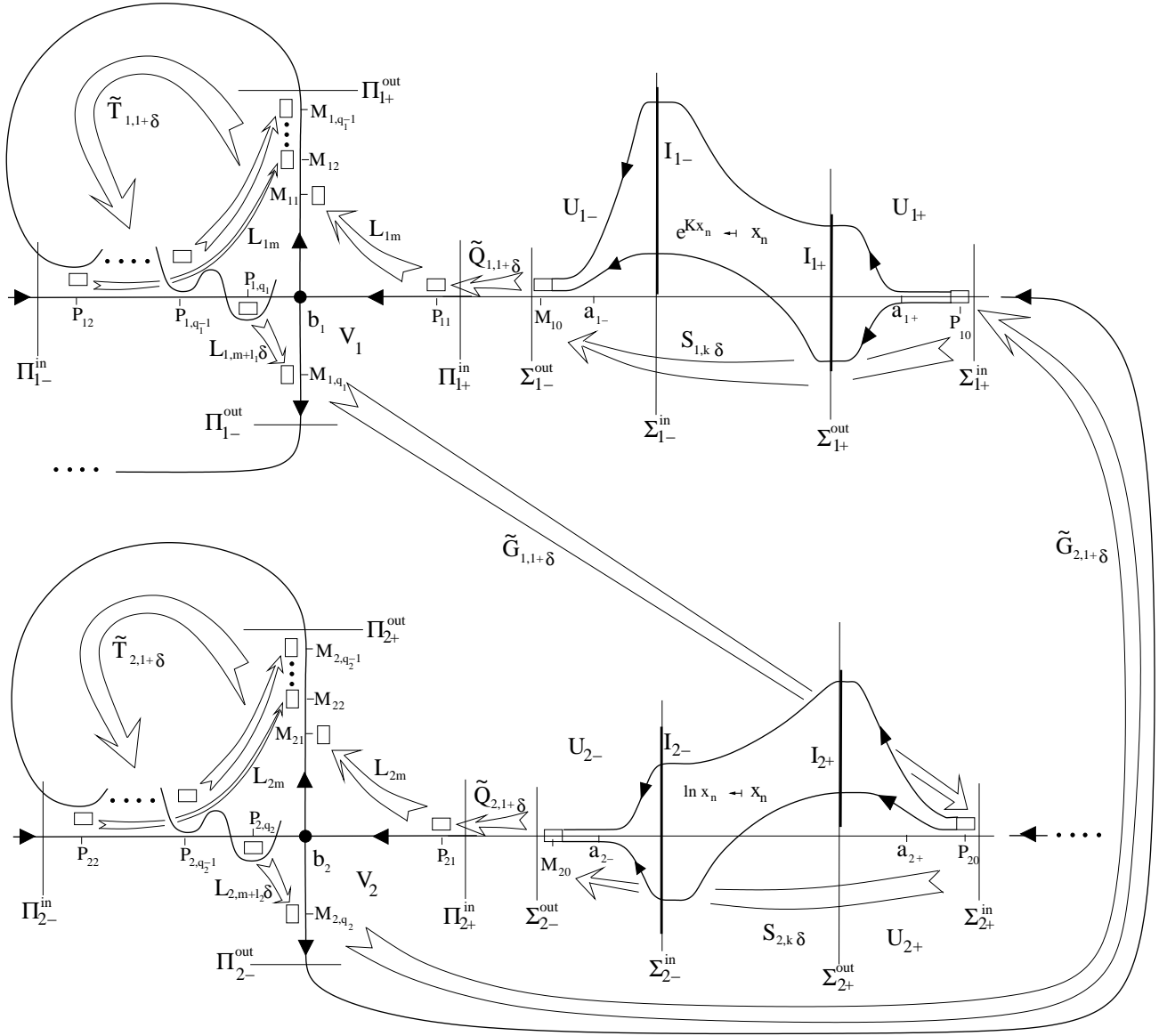}}
\caption{An illustration to the proof of Theorem 1.}
\end{figure}

As $\dot x_1<0$ in $U_{j\sigma}$, every orbit of $Y$ that starts in $U_{j\sigma}$ near
$x_1=a_{j\sigma}+1$ must come in the vicinity of $x_1=a_{j\sigma}-1$ as
time grows. For the corresponding time-$t$ map, we have
\begin{equation}\label{tmap2}
x_1(t)=x_1(0)-t+\frac{1}{2}\beta(\mu_j),\qquad
x_2(t)=e^{\sigma \gamma_j \alpha(\mu_j)} x_2(0),
\end{equation}
where (see (\ref{yw2+}))
\begin{equation}\label{albet}
\begin{array}{l}
\alpha(\mu)=\int_{x_{n-1}(t)-a_{j\sigma}}^{x_{n-1}(0)-a_{j\sigma}}
\frac{\xi(z)dz}{\mu+(1-\mu)(1-\xi(z))}=
\int_{-1/2}^{1/2}\frac{\xi(z)dz}{\mu+(1-\mu)(1-\xi(z))},\\ \\
\beta(\mu)=2\left(\int_{-1/2}^{1/2} \frac{dz}{\mu+(1-\mu)(1-\xi(z))} -1\right).
\end{array}
\end{equation}
Note that $\alpha(\mu)$ is positive, independent of
$x_1(0)$ and $t$ (because we assume that the integration limits
$(x_1(0)-a_{j\sigma})$ and
$(x_1(t)-a_{j\sigma})$ are close to $1$ in the absolute value, i.e.
they fall in the region where $\xi(z)\equiv 0$; see
(\ref{lam})), and both $\alpha(\mu)$ and $\beta(\mu)$
tend to infinity as $\mu\rightarrow+0$ (the integrals diverge at $\mu=0$ because $\xi(0)=1$).

Denote $\Sigma^{in}_{j+}:=\{x_1=a_{j+}+1, |x_2|\leq 1\}$,
$\Sigma^{out}_{j+}:=\{x_1=a_{j+}-1, |x_2|\leq R\}$,
$\Sigma^{in}_{j-}:=\{x_1=a_{j-}+1, |x_2| \leq R\}$,
$\Sigma^{out}_{j-}:=\{x_1=a_{j-}-1, |x_2|\leq 1\}$.
Every orbit of $Y$ that intersects $\Sigma^{in}_{j+}$ at $x_2$ sufficiently small
leaves $U_{j+}$ by crossing $\Sigma^{out}_{j+}$, and the orbits that intersect $\Sigma^{in}_{j-}$
leave $U_{j-}$ by crossing $\Sigma^{out}_{j-}$ (see (\ref{tmap2})).
Define the vector field $Y$ in the region between
$\Sigma^{out}_{j+}$ and $\Sigma^{in}_{j-}$ in such a way that the orbits
starting in $\Sigma^{out}_+$ reach $\Sigma^{in}_-$ at time $1$,
and the corresponding Poincar\'e map
$\Sigma^{out}_{j+}\rightarrow\Sigma^{in}_{j-}$ is
$$x_2\mapsto \psi_j(x_2),$$
where $\psi_1(x_2)\equiv e^{Kx_2}$ at $x_2\in I_{1+}$ and $\psi_2(x_2)\equiv\ln x_2$
at $x_2\in I_{2+}$ (see (\ref{phijs})).
Then, the flow takes the points from the vicinity of $x_1=a_{j+}+1$ in $U_{j+}$ into the vicinity
of $x_1=a_{j-}-1$ in $U_{j-}$. By (\ref{tmap2}), the corresponding time-$t$ map $S_{jt}$ is
\begin{equation}\label{stmap2}
x_1(t)=x_1(0)-t+\beta(\mu_j),\qquad
x_2(t)=e^{-\gamma_j\alpha(\mu_j)}\;
\psi_j\left(e^{\gamma_j\alpha(\mu_j)} x_2(0)\right).
\end{equation}

In the regions $V_j:\{|x_1-b_j|\leq 1, |x_2|\leq 1\}$ ($j=1,2$) we
put the vector field of $Y$ to be equal to
\begin{equation}\label{yw02}
\dot x_1=- (x_1-b_j), \qquad \dot x_2=x_2.
\end{equation}
Thus, in $V_j$, the point $O_j:\{x_1=b_j, x_2=0\}$ is a linear saddle.
Its local stable manifold $W^s_j$ is $x_2=0$,
and the local unstable manifold $W^u_j$ is $x_1=b_j$.
The time-$t$ map $L_{jt}$ within $V_j$ is given by
\begin{equation}\label{ltmap2}
x_1(t)=b_j + e^{-t}(x_1(0)-b_j), \qquad x_2(t)=e^{t} x_2(0).
\end{equation}

Let us define $Y$ in the region between $\Sigma^{out}_{j-}$ and
$\Pi_{j+}^{in}:=\{x_1=b_j+1, |x_2|\leq 1\}$
in such a way that all the orbits starting in a small neighborhood of $x_2=0$ in
$\Sigma^{out}_{j-}$ intersect $\Pi_{j+}^{in}$ at time $1$,
and the corresponding Poincar\'e map is the identity:
$$x_2\mapsto x_2.$$
Then the time-$t$ map $Q_{jt}$ from a small
neighborhood of $x_1=a_{j-}-1, x_2=0$ in $U_{j-}$ into a small neighborhood of
$x_1=b_j+1$ in $V_j$ is given by
\begin{equation}\label{qtmap2}
x_1(t)-b_j=e^{-(t-x_1(0)-2+a_{j-})},\;\; x_2(t)=e^{t-x_1(0)-2+a_{j-}} x_2(0)
\end{equation}
In order to see this, we recall that the vector field in $U_{j-}$ near $x_1=a_{j-}-1$
is given by
$$\dot x_1=-1, \quad \dot x_2=0$$
(see (\ref{yw2+})). Therefore, the term $x_1(0)+2-a_{j-}$ in (\ref{qtmap2})
is the time the orbit spends in order to get from $x(0)$ to $\Pi_{j+}^{in}$).

Every orbit that enters $V_j$ at $x_2>0$ leaves $V_j$ by crossing the cross-section
$\Pi^{out}_{j+}:= \{x_2=1, |x_1-b_j|\leq 1\}$, and
every orbit that enters $V_j$ at $x_2<0$ leaves it by crossing the cross-section
$\Pi^{out}_{j-}:= \{x_2=-1, |x_1-b_j|\leq 1\}$. We assume that
the orbits that start at $\Pi^{out}_{j+}$ close to the point
$W^u_j\cap \Pi^{out}_{j+}=(b_j,1)$ return to $W_j$ at time $1$ and cross
$\Pi^{in}_{j-}:=\{x_1=b_j-1, |x_2|\leq 1\}$; we also assume that
the corresponding Poincar\'e map $x_1\mapsto \bar x_2$ is given by
$$
\bar x_2 = - (x_1-b_j)
$$
(the minus sign stands to ensure the orientability of the flow map). It follows
that the time-$t$ map $T_{jt}$
from a small neighborhood of $W^u_j\cap \Pi^{out}_{j+}$ in $V_{j}$ into a small neighborhood of
$x_1=b_j-1$ in $V_j$ is given by
\begin{equation}\label{ttmap2}
x_1(t)=b_j-e^{-(t-1)}/x_2(0), \qquad x_2(t)=-e^{t-1} x_2(0)^2 (x_1(0)-b_j).
\end{equation}

For the orbits that start at $\Pi^{out}_{j-}$ close to the point
$W^u_j\cap \Pi^{out}_{j-}=(b_j,-1)$ we assume that they cross $\Sigma_{3-j,+}^{in}$ at time $1$,
and the corresponding Poincar\'e map
is given by $\bar x_2 = -(x_1-b_j)$. Thus (see
(\ref{yw2+}),(\ref{lam}),(\ref{yw02})), the time-$t$ map $G_{jt}$
from a small neighborhood of $W^u_j\cap \Pi^{out}_{j-}$
in $V_j$ into a small neighborhood of $x_1=a_{3-j,+}+1$ in $U_{3-j,+}$ is
\begin{equation}\label{gtmap2}
x_1(t)=a_{3-j,+}+2-t-\ln |x_2(0)|,\qquad
x_2(t)=-(x_1(0)-b_j) |x_2(0)|.
\end{equation}

Every $C^\infty$ flow $Y$, which satisfies
(\ref{stmap2}),(\ref{qtmap2}),(\ref{ltmap2}),(\ref{ttmap2}),(\ref{gtmap2}),
is good for our purposes. We may therefore assume that the vector field of $Y$
is identically zero outside some
sufficiently large ball $D$. For small $\delta$, the time-$\delta$ map $Y_\delta$ of the flow
is $O(\delta)$-close to identity in the $C^r$-norm, for any given $r$.
It also equals to identity outside $D$.
Let us fix a certain
$r$, and take a sufficiently small $\delta$ (for convenience,
we assume that $N:=\delta^{-1}$ is an integer). Below we construct an arbitrarily small
(in the $C^r$-norm), localized in $D$ perturbation of $Y_\delta$ as follows.

For the given diffeomorphism $F$, take its sufficiently close approximation in the form of (\ref{dec0n});
in the two-dimensional case the map $\Phi_0$ is given by
\begin{equation}\label{phi02}
\Phi_0:=(x_1,x_2)\mapsto(x_2,-x_1).
\end{equation}
As there is only a finite number ($q_1+q_2$) of the polynomial maps $\tilde H_{js}$ in (\ref{dec0n}),
one can find some finite $d\geq 1$ (common for all $\tilde H_{js}$, $s=1,\dots,q_j, j=1,2$)
such that the maps $\tilde H_{js}$ are written as follows:
\begin{equation}\label{het02}
\bar x_1=x_2, \qquad \bar x_2=-x_1 +
\sum_{0\leq \nu \leq d} h_{js\nu} x_2^{\nu}.
\end{equation}

In the segment $I_j^{out}:=\{e^{-\delta}\leq x_2<1\}$ of $W^u_j$ ($j=1,2$),
we choose $q_j-1$ different points $M_{j1},\dots,M_{j,q_j-1}$,
and one point $M_{jq_j}\in W^u_j$ will be chosen in the segment
$-e^{-\delta}\geq x_2>-1$. Let $u_{js}$
denote the coordinate $x_2$ of $M_{js}$ ($s=1,\dots,q_j$). As $N\delta=1=$flight time
from $\Pi^{out}_{j+}$ to $\Pi^{in}_{j-}$, near the segment $I_j^{out}$
the $(N+1)$-th iteration of the time-$\delta$
map $Y_\delta$ is the map $T_{j,1+\delta}$ from (\ref{ttmap2}), i.e. it is given by
\begin{equation}\label{ydmap2}
\bar x_1=b_j-e^{-\delta}/x_2,\qquad \bar x_2= -e^{\delta} x_2^2 (x_1-b_j).
\end{equation}
This map takes the segment $I_j^{out}$ onto the segment
$\{b_j-1<x_1<b_j-e^{-\delta}, x_2=0\}\in W^s_j$. Let
$P_{j,s+1}=T_{j,1+\delta}M_{js}$
($s=1,\dots,q_j-1$), and let $P_{j1}$ be a point from
$\{b_j+e^{-\delta}<x_1<b_j+1, x_2=0\}\in W^s_j$.
We denote the coordinate $x_1$ of $P_{j,s+1}$ as $z_{js}$.
By (\ref{ydmap2}),
\begin{equation}\label{auxz2}
z_{j,s+1}=b_j-e^{-\delta}/u_{js}.
\end{equation}
We also assume
\begin{equation}\label{z12}
z_{j1}=b_j+e^{-\delta/2}.
\end{equation}
We stress that $u_{js}$ and $(b_j-z_{js})$ are bounded away from zero.

Take a sufficiently large integer $m$ and choose some points $P_{js}'=(z_{js},z_{js}')$ and
$M_{js}'=(u_{js}',u_{js})$,
sufficiently close to $P_{js}$ and $M_{js}$ respectively ($j=1,2; s=1,\dots, q_j$).
We define the coordinates of $P_{js}'$ and $M_{js}'$ by the following rule:
\begin{equation}\label{pmuz2}
\begin{array}{l}\displaystyle
z_{js}'= e^{-m} u_{js}),\quad u_{js}'=b_j+e^{-m}(z_{js}-b_j)\;\;
\mbox{ at }\; s\leq q_j-1\\ \displaystyle
z_{jq_j}'= e^{-(m+l_j\delta)} u_{jq_j},\quad
u_{jq_j}'=b_j+e^{-(m+l_j\delta)}(z_{jq_j}-b_j),
\end{array}
\end{equation}
where $l_j$ is an integer to be defined later
(see (\ref{ljdf2});
note that $l_j\delta$ remains uniformly bounded for arbitrarily small $\delta$).
As $m$ is assumed to be large and $l_j\delta$ is bounded, such defined points $P_{js}'$
and $M_{js}'$ are closed to respective points $P_{js}$ and $M_{js}$ indeed.
It follows from (\ref{pmuz2}) that $M_{js}'=L_{jm} P_{js}'$ at $s\leq q_j-1$,
where $L_{jt}$ is the map (\ref{ltmap2}). At $s=q_j$ we have
$M_{jq_j}'=L_{j,m+l_j\delta} P_{jq_j}'$.

Let us add to the map $Y_\delta$ a small perturbation, which is
localized in a small neighborhood of the points
$Y_\delta^{-1}(P_{j,s+1})$ (so outside these small neighborhoods $Y_\delta$ remains unchanged).
We require that these localized perturbations are such that in a sufficiently
small neighborhood of $M_{js}=Y_{\delta}^{-N}(Y_\delta^{-1}P_{j,s+1})$ the map
$\tilde Y_{\delta}^{N+1}$
(where $\tilde Y_{\delta}$ denotes the perturbed map) is given by the following
perturbation of (\ref{ydmap2}):
\begin{equation}\label{pcr2}
\begin{array}{l}\displaystyle
\bar x_1=b_j-e^{-\delta}/x_2,\\ \displaystyle
\bar x_2= z_{j,s+1}'-e^{\delta} x_2^2 (x_1-u_{js}') + \;\sum_{0\leq \nu \leq d}
\varepsilon_{js\nu} (x_2-u_{js})^{\nu},\end{array}
\end{equation}
where $\varepsilon_{js\nu}$ are small coefficients
to be determined later (see (\ref{epj2})). By (\ref{pmuz2}), $z_{j,s+1}'$ and $(u_{js}'-b_j)$
are small as well, so (\ref{pcr2}) is a small perturbation of (\ref{ydmap2}) indeed.
It is easy to see from (\ref{pcr2}),(\ref{auxz2}) that
$P_{j,s+1}'=\tilde Y_\delta^{N+1} M_{js}'$  at $\varepsilon=0$.

At $\varepsilon=0$ (hence at all small $\varepsilon$) the map
$\tilde T_{j,1+\delta}\circ L_{jm}\equiv \tilde Y_{\delta}^{mN+N+1}$ (where $\tilde T$ stands
for the perturbed map $T$) takes a small neighborhood of $P_{js}'=(z_{js},z_{js}')$ into a small
neighborhood of $P_{j,s+1}'=(z_{j,s+1},z_{j,s+1}')$. By (\ref{pcr2}),(\ref{ltmap2}),(\ref{pmuz2}),
this map is written as
\begin{equation}\label{tl2m}
\begin{array}{l}\displaystyle
\bar x_1=b_j-e^{-(m+\delta)}/x_2,\\ \displaystyle
\bar x_2= z_{j,s+1}'-e^{(m+\delta)} x_2^2 (x_1-z_{js}) + \;\sum_{0\leq \nu \leq d}
\varepsilon_{js\nu} e^{\nu m} (x_2-z_{js}')^{\nu},\end{array}
\end{equation}

We choose some $\eta(m)$ that tends to zero as $m\rightarrow+\infty$ and
introduce rescaled coordinates $(v_1,v_2)$ near $P_{js}'$ by the rule
\begin{equation}\label{scal2}
x_1 = z_{js} + C_{js}\; \eta \; v_1,\qquad
x_2 = z_{js}'+ \frac{1}{C_{js}} \; \eta \; e^{-m} v_2,
\end{equation}
where the independent of $m$ positive coefficients $C_{js}$ are determined
later (see (\ref{cdef2}); note that $C_{js}$ are bounded away from zero and infinity).
Since $\eta$ tends to zero as $m\rightarrow+\infty$, any bounded region of values
of $v$ corresponds to a small neighborhood of $P_{js}'$.

After the rescaling, map (\ref{tl2m}) takes the following form
(see (\ref{pmuz2}),(\ref{auxz2})):
$$\begin{array}{l} \displaystyle
C_{j,s+1} \bar v_1= \frac{1}{e^{\delta}u_{js}^2 C_{js}} v_2+O(\eta),\\
\displaystyle
\frac{1}{C_{j,s+1}} \bar v_2 = - e^{\delta} u_{js}^2 C_{js} v_1
+ O(\eta)+ \;\sum_{0\leq \nu \leq d}
\varepsilon_{js\nu} e^{m}\eta^{\nu-1} C_{js2}^{\nu} v_2^{\nu},
\end{array}
$$
As we see, by putting
\begin{equation}\label{cdef2}
C_{j,s+1}=\frac{1}{e^{\delta}u_{js}^2 C_{js}},
\end{equation}
and
\begin{equation}\label{epj2}
\varepsilon_{js\nu}=h_{js\nu} e^{-m} \eta^{1-\nu} C_{j,s+1}^{-1} C_{js}^{\nu},
\end{equation}
the map $\tilde T_{j,1+\delta}\circ L_{jm}$ near $P_{js}'$ takes the form
\begin{equation}\label{scsdl2}
\begin{array}{l}
\bar v_1= v_2 + O(\eta),\\
\bar v_2 = -v_1 + O(\eta)
+\;\sum_{0\leq\nu \leq d}
h_{js\nu} v_2^{\nu},
\end{array}
\end{equation}
i.e. it can be made as close as we want to the map $\tilde H_{js}$, provided $m$
is taken large enough (recall that $\eta\rightarrow 0$ as $m\rightarrow+\infty$).
We take $\eta$ tending to zero sufficiently
slowly, so all $\varepsilon_{js\nu}\rightarrow 0$ (see (\ref{epj2}); recall that $\nu\leq d$
where $d$ is independent of $m$ and $\delta$). Thus,
our perturbation to $Y_\delta$ is arbitrarily small indeed.

It follows that in the rescaled coordinates the map
$\left(\tilde T_{j,1+\delta}\circ L_{jm}\right)^{q_j}\equiv \tilde Y_\delta^{q_j(mN+N+1)}$
from a small neighborhood of $P_{j1}'$ into a small neighborhood of $P_{jq_j}'$
can be made as close as we want to the map
$$\tilde H_{jq_j}\circ \dots \circ \tilde H_{j1},$$
provided $m$ is large enough. The rescaled coordinates near $P_{j1}'$ and $P_{jq_j}'$ are
given by formulas (\ref{scal2}), where the coefficients $C_{ji}>0$ may be taken arbitrary,
and the coefficients $C_{jq_j}$ are then recovered from the recursive formula (\ref{cdef2}).
Further it is convenient to put
\begin{equation}\label{c1fx}
C_{j1}=e^{-\delta/2}.
\end{equation}
Note that $m$ does not enter (\ref{cdef2}),(\ref{c1fx}),
hence $C_{jq_j}$ stay bounded away from zero and infinity as $m\rightarrow+\infty$.

Let us also fix the choice of the integer $l_j$ that enters the definition of the points
$P_{jq_j}'$ and $M_{jq_j}'$ (see (\ref{pmuz2})). Namely, we require that
\begin{equation}\label{ljdf2}
\ln|u_{jq_j}|+\ln C_{jq_j} = l_j\delta.
\end{equation}
Recall that $u_{jqj}$ is the coordinate $x_2$ of the point $M_{jq_j}$, and
it can be arbitrarily taken within the interval
\begin{equation}\label{ujqjb}
-e^{-\delta}\geq u_{jqj}>-1.
\end{equation}
As $C_{jq_j}$ does not depend on the choice of $u_{jqj}$ (see (\ref{cdef2})),
condition (\ref{ljdf2}) uniquely defines both the integer $l_j$ and the value of $u_{jq_j}$.
The claimed before uniform boundedness of $l_j\delta$ follows from the uniform
boundedness of $\ln C_{jq_j}$.

Now, from (\ref{ltmap2}),(\ref{scal2}),(\ref{pmuz2}),(\ref{ljdf2}) we obtain that the map
$L_{j,m+l_j\delta}$
from a small neighborhood of $P_{jq_j}'$ into a small neighborhood of $M_{jq_j}'$
is identity (i.e. $\bar v_1=v_1,\; \bar v_2 =v_2$),
provided the rescaled coordinates $(v_1,v_2)$ near $M_{jq_j}'$ are introduced as follows:
\begin{equation}\label{scall2}
\begin{array}{l}\displaystyle
x_1=u_{jq_j}'+|u_{jq_j}|^{-1} \; \eta \; e^{-m} v_1,\\
\displaystyle x_2=u_{jq_j} - u_{jq_j} \; \eta \;  v_2.
\end{array}
\end{equation}
Therefore, for the given choice of the coordinates, by taking $m\rightarrow+\infty$, the map
$L_{j,m+l_j\delta}\circ\left(\tilde T_{j,1+\delta}\circ L_{jm}\right)^{q_j}\equiv
\tilde Y_\delta^{l_j+mN+q_j(mN+N+1)}$
from a small neighborhood of $P_{j1}'$ into a small neighborhood of $M_{jq_j}'$
can be made as close as we want to $\tilde H_{jq_j}\circ \dots \circ \tilde H_{j1}$
(we have already proved the same for the map
$\left(\tilde T_{j,1+\delta}\circ L_{jm}\right)^{q_j}$, and the map
$L_{j,m+l_j\delta}$ is identity in the chosen coordinates).

Recall that, by construction, the point $Y_\delta^{N+1}M_{jq_j}'$ lies
in $U_{3-j,+}$ in the region $a_{3-j,+}+1-\delta<x_1<a_{3-j,+}+1$. We add to
the map $Y_\delta$ a perturbation, localized near the point $Y_\delta^{N}M_{jq_j}'$, such that
the corresponding map $\tilde G_{j,1+\delta}\equiv \tilde Y_\delta^{N+1}$ will
have the following form near $M_{jq_j}'$:
\begin{equation}\label{gtmapt2}
\bar x_1=a_{3-j,+}+1-\delta-\ln |x_2|,\qquad
\bar x_2=-|x_2| (x_1-u_{jq_j}');
\end{equation}
since $u_{jq_j}'\rightarrow b_j$ as $m\rightarrow+\infty$ (see (\ref{pmuz2})), these
fomulas define a small perturbation of the map
$G_{j,1+\delta}$ given by (\ref{gtmap2}) indeed.

Denote $P_{3-j,0}'=\tilde G_{j,1+\delta} M_{jq_j}'$. By (\ref{gtmapt2}), the
coordinates of $P_{3-j,0}'$ are given by $\{x_1=a_{3-j,+}+1-\delta-\ln|u_{jq_j}|,\; x_2=0\}$.
Introduce rescaled coordinates near the points $P_{3-j,0}'$ ($j=1,2$) by the rule
\begin{equation}\label{resp32}
x_1=a_{3-j,+}+1-\delta-\ln|u_{jq_j}|+ \eta \; v_1,\qquad x_2=\eta \; e^{-m} v_2.
\end{equation}
When coordinates are rescaled by rule (\ref{scall2}) near $M_{jq_j}'$ and by rule (\ref{resp32})
near $P_{3-j,0}'$, map (\ref{gtmapt2}) takes the form
$(v_1,v_2)\mapsto(v_2,-v_1)+O(\eta)$, i.e. it becomes arbitrarily
close to the map $\Phi_0$ (see (\ref{phi02})) as $m_\rightarrow+\infty$. Thus, in the
coordinates rescaled by rule (\ref{scal2}) near $P_{j1}'$ and by rule (\ref{resp32})
near $P_{3-j,0}'$, the map
$\tilde G_{j,1+\delta}\circ L_{j,m+l_j\delta}\circ
\left(\tilde T_{j,1+\delta}\circ L_{jm}\right)^{q_j}\equiv
\tilde Y_\delta^{l_j+(1+q_j)(mN+N+1)}$
from a small neighborhood of $P_{j1}'$ into a small neighborhood of $P_{3-j,0}'$ ,
can be made as close as we want to the map
$\Phi_0\circ\tilde H_{jq_j}\circ \dots \circ \tilde H_{j1}$
as $m$ grows.

Analogously, we take the point $M_{j0}':\{x_1=a_{j-}-1+\delta/2,\; x_2=0 \}\in U_{j-}$
and perturb the map $Y_\delta$ near $Y_\delta^N P_{j0}'$ in such a way that the map
$\tilde Q_{j,1+\delta}\equiv \tilde Y_\delta^{N+1}$ near $M_{j0}'$ will be given by
\begin{equation}\label{qtmapm2}
\begin{array}{l} \displaystyle
\bar x_1-b_j=e^{-(\delta-x_1-1+a_{j-})},\;\;
\bar x_2=e^{\delta-x_1-1+a_{j-}} x_2 + e^{-m}u_{j1}.
\end{array}
\end{equation}
It is a small perturbation of the map $Q_{j,1+\delta}$ given by (\ref{qtmap2}), and it takes
$M_{j0}'$ to $P_{j1}'$ (see (\ref{z12}),(\ref{pmuz2})). When we introduce rescaled variables
near $M_{j0}'$ by the rule
\begin{equation}\label{resp02}
x_1=a_{j-}-1+\delta/2+ \eta\; v_1,\qquad x_2= \eta \; e^{-m} v_2,
\end{equation}
and near $P_{j1}'$ by rule (\ref{scal2}),(\ref{c1fx}), map (\ref{qtmapm2}) will take the
form $\bar v=v+O(\eta)$, i.e. it is close to the identity map.
Thus, in the rescaled coordinates given by (\ref{resp02}),(\ref{resp32}), the map
$\tilde G_{j,1+\delta}\circ L_{j,m+l_j\delta}\circ
\left(\tilde T_{j,1+\delta}\circ L_{jm}\right)^{q_j}\circ\tilde Q_{j,1+\delta} \equiv
\tilde Y_\delta^{l_j+(1+q_j)(mN+N+1)+N+1}$
from a small neighborhood of $M_{j0}'$ into a small neighborhood of $P_{3-j,0}'$,
is as close as we want to the map $\Phi_0\circ\tilde H_{jq_j}\circ \dots \circ \tilde H_{j1}$
at $m$ large enough, i.e. it is a close approximation of the map $\Phi_j$.

Let us now determine the form of the map $S_{jt}:v\mapsto\bar v$
from a small neighborhood of $P_{j0}'$ into a small neighborhood of $M_{j0}'$
in the coordinates rescaled by the same rules (\ref{resp32}) and (\ref{resp02})
(note that when rescaling near $P_{j0}'$ one should change $j$ to $3-j$
in the right-hand side of (\ref{resp32})). We will
choose the coefficients $\mu_j$ in (\ref{yw2+}) such that
\begin{equation}\label{mdf2}
\beta(\mu_j)= k-5+3/2\delta+\ln|u_{3-j,q_{3-j}}|
\end{equation}
where $k$ is some sufficiently large integer. Since
$\beta\rightarrow+\infty$ as $\mu\rightarrow+0$ (see (\ref{albet})),
and $\ln|u_{3-j,q_{3-j}}|$ is bounded (see (\ref{ujqjb})),
equation (\ref{mdf2}) has a solution $\mu_j(k)$ for every sufficiently large $k$,
and $\mu_j(k)\rightarrow +0$ as $k\rightarrow+\infty$.
It follows that
$\alpha(\mu_j(k))\rightarrow+\infty$. Thus, for any sufficiently
large $m$ we can find $\gamma_j\in (0,1]$ and large $k$ such that
\begin{equation}\label{mdfa2}
\eta(m) \; e^{-m}=e^{-\gamma_j \alpha(\mu_j(k))}.
\end{equation}
From (\ref{stmap2}) we immediately obtain
that with this choice of $\mu_j$, $\gamma_j$ and $k$
the map $S_{j,k\delta}\equiv Y_{\delta}^{k}$ from a small neighborhood
the point $P_{j0}'$ into  a small neighborhood of $M_{j0}'$
takes the following form in the coordinates rescaled by rules, respectively,
(\ref{resp32}) (with $(3-j)$ changed to $j$) and (\ref{resp02}):
$$\bar v_1=v_1, \qquad \bar v_2=\psi_j(v_2).$$
As wee see, the map $S_{j,k\delta}$ in the rescaled coordinates coincides with the map $\Psi_j$
for $v$ from some open neighborhood of $D$ (if $j=1$) or of $\Phi_1\circ\Psi_1(D)$ (if $j=2$).

Summarizing, we obtain that the map
$$\begin{array}{l} \displaystyle
\tilde G_{2,1+\delta}\circ L_{2,m+l_2\delta}\circ
\left(\tilde T_{2,1+\delta}\circ L_{2m}\right)^{q_2}\circ\tilde Q_{2,1+\delta}\circ S_{2,k\delta}\circ\\
\displaystyle
\qquad\qquad \circ\; \tilde G_{1,1+\delta}\circ L_{1,m+l_1\delta}\circ
\left(\tilde T_{1,1+\delta}\circ L_{1m}\right)^{q_1}\circ\tilde Q_{1,1+\delta}\circ S_{1,k\delta}\equiv\\
\displaystyle \equiv \;
\tilde Y_\delta^{2k+l_1+l_2+(2+q_1+q_2)(mN+N+1)+2(N+1)}
\end{array}$$
is a close approximation to the map $\Phi_2\circ\Psi_2\circ\Phi_1\circ\Psi_1$
(i.e. to the original map $F$),
provided $\tilde H_{jq_j}\circ\dots \circ \tilde H_{j1}$ are sufficiently close approximations to
$\tilde \Phi_j$ ($j=1,2$) and $m$ is large enough. This completes the proof of the theorem
in the two-dimensional case.

\subsection{Higher-dimensional case}
\setcounter{equation}{0}

In the case of arbitrary $n>2$, the construction follows the same line as in the two-dimensional case.
As above, $\Phi_{1,2}$ and $\Psi_{1,2}$ are the maps defined by (\ref{decom}), and $I_{1\pm}$ and
$I_{2\pm}$ are intervals of values of $x_n$ such that $x_n\in I_{1+}$ at $x\in B^n$,
$x_n\in I_{1-}$ at $x\in \Psi_1(B^n)$, $x_n\in I_{2+}$ at $x\in \Phi_1\circ\Psi_1(B^n)$
and $x_n\in I_{2-}$ at $x\in \Psi_2\circ\Phi_1\circ\Psi_1(B^n)$. A value $R$ is chosen such that all the
intervals $I_{j\pm}$ lie within $\{|x_n|\leq R\}$.
Choose numbers $a_{1+}=a_{1-}+3=b_1+6=a_{2+}+9=a_{2-}+12=b_2+15$. Define the regions
$U_{j\sigma}:\{|x_{n-1}-a_{j\sigma}\|\leq 1, |x_n|\leq R, |x_i|\leq 1 \; (i\leq n-2)\}$,
and $V_j:\{|x_{n-1}-b_j|\leq 1, |x_i|\leq 1 \;\; (i\neq n-1)\}$, $j=1,2$, $\sigma=\pm 1$.
Let the vector field of a $C^\infty$ flow $Y$ in $U_{j\sigma}$ be equal to
\begin{equation}\label{yw+}
\begin{array}{l}\displaystyle
\dot x_{n-1}=-\mu_j-(1-\mu_j)(1-\xi(x_{n-1}-a_{j\sigma})),\\ \\\displaystyle
\dot x_i=\sigma \gamma_{i\sigma} \; x_i \;\xi(x_{n-1}-a_{j\sigma}) \;\;\; (i\neq n-1),
\end{array}
\end{equation}
where $\mu_{1,2}>0$ are small (see (\ref{mdf}),(\ref{mdfa})), $\gamma_{i\pm}\in[0,1]$ (see
(\ref{mdfa}),(\ref{gdf})), and the $C^\infty$ function $\xi$ satisfies (\ref{lam}).
In the regions $V_j$ we make $Y$ equal to
\begin{equation}\label{yw0}
\dot x_i=-\lambda_i x_i \;\; (i=1,\dots,n-2), \;\;\dot x_{n-1}=-\lambda_{n-1} (x_{n-1}-b_j),
\quad \dot x_n=x_n;
\end{equation}
here $\lambda_i>0$ are such that
\begin{equation}\label{gam}
\lambda_2=\dots=\lambda_{n-1}=\lambda, \qquad \lambda_1=1-(n-2)\lambda,
\end{equation}
where the positive number $\lambda$ is specified below (see (\ref{gamd})).

As $\dot x_{n-1}<0$ in $U_{j\sigma}$, every orbit of $Y$ that starts in $U_{j\sigma}$ near
$x_{n-1}=a_{j\sigma}+1$ must come in the vicinity of $x_{n-1}=a_{j\sigma}-1$ as
time grows. For the corresponding time-$t$ map, we have
\begin{equation}\label{tmap}
x_i(t)=e^{\sigma \gamma_{i\sigma} \alpha(\mu_j)} x_i(0) \;\;\; (i\leq n-2), \qquad
x_n(t)=x_n(0)-t+\frac{1}{2}\beta(\mu_j),
\end{equation}
where the tending to infinity, as $\mu\rightarrow+0$, functions $\alpha(\mu)$ and $\beta(\mu)$
are defined by (\ref{albet}).

Denote $\Sigma^{in}_{j+}:=\{x_{n-1}=a_{j+}+1, |x_n|\leq 1\}$,
$\Sigma^{out}_{j+}:=\{x_{n-1}=a_{j+}-1, |x_n|\leq R\}$,
$\Sigma^{in}_{j-}:=\{x_{n-1}=a_{j-}+1, |x_n| \leq R\}$,
$\Sigma^{out}_{j-}:=\{x_{n-1}=a_{j-}-1, |x_n|\leq 1\}$
(we also assume that $|x_i|\leq 1$ for $i\leq n-2$ on $\Sigma_{j\pm}^{in,out}$).
Every orbit of $Y$ that intersects $\Sigma^{in}_{j+}$ at $x_n,x_1,\dots,x_{n-2}$ sufficiently small
leaves $U_{j+}$ by crossing $\Sigma^{out}_{j+}$, and the orbits that intersect $\Sigma^{in}_{j-}$
leave $U_{j-}$ by crossing $\Sigma^{out}_{j-}$ (see (\ref{tmap})).
We define $Y$ in the region between
$\Sigma^{out}_{j+}$ and $\Sigma^{in}_{j-}$ in such a way that the orbits starting in $\Sigma^{out}_+$
reach $\Sigma^{in}_-$ at time $1$,
and the corresponding Poincar\'e map
$\Sigma^{out}_{j+}\rightarrow\Sigma^{in}_{j-}$ is
$(x_1,\dots,x_{n-2},x_n)\mapsto(x_1,\dots,x_{n-2},\psi_j(x_n))$,
where the functions $\psi_j$ are defined by (\ref{phijs}).
Then, the flow takes the points from the vicinity of $x_{n-1}=a_{j+}+1$ in $U_{j+}$ into the vicinity
of $x_{n-1}=a_{j-}-1$ in $U_{j-}$. By (\ref{tmap}), the corresponding time-$t$ map $S_{jt}$ is
\begin{equation}\label{stmap}
\begin{array}{l}
x_i(t)=e^{(\gamma_{_{i+}}-\gamma_{_{i-}})\alpha(\mu_j)} x_i(0) \qquad (i\leq  n-2),\\
x_{n-1}(t)=x_{n-1}(0)-t+\beta(\mu_j),\qquad
x_n(t)=e^{-\gamma_{_{n-}}\alpha(\mu_j)}\psi_j(e^{\gamma_{_{n+}}\alpha(\mu_j)} x_n(0)).
\end{array}
\end{equation}

In the region between $\Sigma^{out}_{j-}$ and $\Pi_{j+}^{in}:=\{x_{n-1}=b_j+1, |x_i|\leq 1\}$, we define $Y$
in such a way that all the orbits starting in a small neighborhood of $x_n=x_1=\dots=x_{n-2}=0$ in
$\Sigma^{out}_{j-}$ intersect $\Pi_{j+}^{in}$ at time $1$, and the corresponding Poincar\'e map
is the identity:
$(x_1,\dots,x_{n-2},x_n)\mapsto(x_1,\dots,x_{n-2},x_n)$. Then the time-$t$ map $Q_{jt}$ from a small
neighborhood of $x_{n-1}=a_{j-}-1, x_n=x_1=\dots=x_{n-2}=0$ in $U_{j-}$ into a small neighborhood of
$x_{n-1}=b_j+1$ in $V_j$ is given by
\begin{equation}\label{qtmap}
\begin{array}{l}
x_i(t)=e^{-\lambda_i(t-x_{n-1}(0)-2+a_{j-})}x_i(0) \qquad (i\leq  n-2),\\
x_{n-1}(t)-b_j=e^{-\lambda(t-x_{n-1}(0)-2+a_{j-})},\;\;
x_n(t)=e^{t-x_{n-1}(0)-2+a_{j-}} x_n(0)
\end{array}
\end{equation}
(see (\ref{yw0})).

In $V_j$ ($j=1,2$), the local stable manifold $W^s_j$ of the linear saddle equilibrium state
$O_j:\{x_{n-1}=b_j, x_i=0 \; (i\neq n-1)\}$ is $x_n=0$, and the local unstable manifold $W^u_j$
is $x_{n-1}=b_j, x_1=\dots=x_{n-2}=0$. The time-$t$ map $L_{jt}$ within $V_j$ is given by
\begin{equation}\label{ltmap}
\begin{array}{l}
x_i(t)=e^{-\lambda_i t}x_i(0) \qquad (i\leq  n-2),\\
 x_{n-1}(t)-b_j=e^{-\lambda t}(x_{n-1}(0)-b_j), \qquad x_n(t)=e^{t} x_n(0).
 \end{array}
\end{equation}
Every orbit that enters $V_j$ at $x_n>0$ leaves $V_j$ by crossing the cross-section
$\Pi^{out}_{j+}:= \{x_n=1, |x_{n-1}-b_j|\leq 1, |x_i|\leq 1 \; (i \leq n-2)\}$, and
every orbit that enters $V_j$ at $x_n<0$ leaves it by crossing the cross-section
$\Pi^{out}_{j-}:= \{x_n=-1, |x_{n-1}-b_j|\leq 1, |x_i|\leq 1 \; (i \leq n-2)\}$. We assume that
the orbits that start at $\Pi^{out}_{j+}$ close to the point
$W^u_j\cap \Pi^{out}_{j+}=(x_1=\dots=x_{n-2}=0, x_{n-1}=b_j)$ return to $W_j$ at time $1$ and cross
$\Pi^{in}_{j-}:=\{x_{n-1}=b_j-1, |x_i|\leq 1 (i\neq n-1)\}$; we also assume that
the corresponding Poincar\'e map $(x_1,\dots,x_{n-1})\mapsto(\bar x_1,\dots,\bar x_{n-2},\bar x_n)$
is given by
$$
\bar x_i = x_{i+1} \;\;\; (i\leq n-3),\qquad \bar x_{n-2}=x_{n-1}-b_j, \qquad \bar x_n = (-1)^{n+1} x_1
$$
(the factor $(-1)^{n+1}$ stands to ensure the orientability). It follows that the time-$t$ map $T_{jt}$
from a small neighborhood of $W^u_j\cap \Pi^{out}_{j+}$ in $V_{j}$ into a small neighborhood of
$x_{n-1}=b_j-1$ in $V_j$ is given by
\begin{equation}\label{ttmap}
\begin{array}{l}
x_i(t)=e^{-\lambda_i (t-1)} x_n(0)^{\lambda_{i+1}-\lambda_i} x_{i+1}(0) \qquad (i\leq  n-3),\\
x_{n-2}(t)=e^{-\lambda_{n-2} (t-1)} x_n(0)^{\lambda_{n-1}-\lambda_{n-2}} (x_{n-1}(0)-b_j),\\
x_{n-1}(t)=b_j-e^{-\lambda_{n-1}(t-1)} x_{n}(0)^{-\lambda_{n-1}}, \\
x_n(t)=(-1)^{n+1} e^{t-1} x_{n}(0)^{1+\lambda_1} x_1(0).
\end{array}
\end{equation}

For the orbits that leave $V_j$ by crossing $\Pi^{out}_{j-}$, we assume that
the orbits that start at $\Pi^{out}_{j-}$ close to the point
$W^u_j\cap \Pi^{out}_{j-}=(x_1=\dots=x_{n-2}=0, x_{n-1}=b_j)$ cross $\Sigma_{3-j,+}^{in}$ at time $1$,
and the corresponding Poincar\'e map $(x_1,\dots,x_{n-1})\mapsto(\bar x_1,\dots,\bar x_{n-2},\bar x_n)$
is given by $\bar x_i = x_{i}$ at $i=1,\dots,n-2$ and $\bar x_n = -(x_{n-1}-b_j)$. Thus (see
(\ref{yw+}),(\ref{lam}),(\ref{yw0})), the time-$t$ map $G_{jt}$
from a small neighborhood of $W^u_j\cap \Pi^{out}_{j-}$
in $V_j$ into a small neighborhood of $x_{n-1}=a_{3-j,+}+1$ in $U_{3-j,+}$ is
\begin{equation}\label{gtmap}
\begin{array}{l}
x_i(t)=|x_n(0)|^{\lambda_i} x_i(0) \qquad (i\leq  n-2),\\
x_{n-1}(t)=a_{3-j,+}+2-t-\ln |x_n(0)|,\\
x_n(t)=-|x_{n}(0)|^{\lambda_{n-1}} (x_{n-1}(0)-b_j).
\end{array}
\end{equation}

These conditions define the flow $Y$ (we also assume that the vector field of $Y$ is identically
zero outside some sufficiently large ball $D$). Let us take a sufficiently small $\delta$ such
that $N:=\delta^{-1}$ is an integer, and proceed to the construction of a small
(arbitrarily small in the $C^r$-norm, with any chosen in advance $r$) perturbation of $Y_\delta$,
localized in $D$.

Take a sufficiently close approximation of the given diffeomorphism $F$ by the product (\ref{dec0n}).
There exists some finite $d\geq 1$, common for all $\tilde H_{js}$ ($s=1,\dots,q_j, j=1,2$)
such that the polynomial maps $\tilde H_{js}$ in (\ref{dec0n}) are written as follows:
\begin{equation}\label{het0}
\bar x_i=x_{i+1} \;\;\; (i\leq n-1), \qquad \bar x_n=(-1)^{n+1}x_1 +
\sum_{\stackrel{\nu_2\geq 0,\dots,\nu_n\geq 0}{
\nu_2+\dots+\nu_n \leq d}}
h_{js\nu} \prod_{2\leq p\leq n} x_p^{\nu_p}.
\end{equation}
We will now fix the choice of the positive $\lambda$ in (\ref{gam}) such that
\begin{equation}\label{gamd}
\lambda<\frac{1}{(n-1) d + r}.
\end{equation}

In the segment $I_j^{out}:=\{e^{-\delta}\leq x_n<1\}$ of $W^u_j$ ($j=1,2$),
we choose $q_j-1$ different points $M_{j1},\dots,M_{j,q_j-1}$,
and one point $M_{jq_j}\in W^u_j$ will be chosen in the segment $-e^{-\delta}\geq x_n>-1$. Let $u_{js}$
denote the coordinate $x_n$ of $M_{js}$ ($s=1,\dots,q_j$). As $N\delta=1$,
near the segment $I_j^{out}$ the $(N+1)$-th iteration of the time-$\delta$
map $Y_\delta$ is the map $T_{j,1+\delta}$ from (\ref{ttmap}).
Thus, the map $Y_\delta^{N+1}$ near $I_j^{out}$ will be given by
\begin{equation}\label{ydmap}
\begin{array}{l} \displaystyle
\bar x_1=e^{((n-2)\lambda-1)\delta} x_n^{(n-1)\lambda-1} \hat x_{2},\qquad
\bar x_i=e^{-\lambda \delta} \hat x_{i+1} \;\;\; (2\leq i\leq  n-2),\\ \\ \displaystyle
\bar x_{n-1}=b_j-e^{-\lambda\delta} x_{n}^{-\lambda},\qquad
\bar x_n= (-1)^{n+1} e^{\delta} x_{n}^{2-(n-2)\lambda} x_1
\end{array}
\end{equation}
(where we denote $\hat x_i=x_i$ at $i\neq n-1$ and $\hat x_{n-1}=x_{n-1}-b_j$).
This map takes the segment $I_j^{out}$ onto the segment
$\{b_j-1<x_{n-1}<b_j-e^{-\lambda\delta}, x_1=\dots=x_{n-2}=x_n=0\}\in W^s_j$. Let
$P_{j,s+1}=T_{j,1+\delta}M_{js}$
($s=1,\dots,q_j-1$), and let $P_{j1}$ be a point from
$\{b_j+e^{-\lambda\delta}<x_{n-1}<b_j+1,x_1=\dots=x_{n-2}=x_n=0\}\in W^s_j$.
By (\ref{ydmap}), the coordinate $x_{n-1}$ of $P_{j,s+1}$ equals to
$b_j-e^{-\lambda\delta} u_{js}^{-\lambda}$.

We take sufficiently large integer $m$ and choose some points $P_{js}'$ and $M_{js}'$, sufficiently close
to $P_{js}$ and $M_{js}$ respectively ($j=1,2; s=1,\dots, q_j$), such that at $s\leq q_j-1$ we have
$M_{js}'=L_{jm} P_{js}'$ (where $L_{jt}$ is the map (\ref{ltmap})). At $s=q_j$ we assume
$M_{jq_j}'=L_{j,m+l_j\delta} P_{jq_j}'$ where $l_j$ is an integer to be defined later (see (\ref{ljdf});
note that $l_j\delta$ is uniformly bounded). Denote the coordinates
of $P_{js}'$ and $M_{js}'$ as $(z_{js1}',\dots,z_{js,n-2}',b_j+z_{js,n-1}',z_{jsn}')$ and
$(u_{js1}',\dots,u_{js,n-2}',b_j+u_{js,n-1}',u_{jsn}')$ respectively. By (\ref{ltmap}),
\begin{equation}\label{uzr}
u_{jsi}'=e^{-\lambda_i m} z_{jsi}' \;\; (i=1,\dots,n-1), \;\; u_{jsn}'=e^{m} z_{jsn}'
\end{equation}
at $s\leq q_j-1$. At $s=q_j$ we have
\begin{equation}\label{uzrl}
u_{jq_ji}'=e^{-\lambda_i (m+l_j\delta)} z_{jq_ji}' \;\;\; (i=1,\dots,n-1), \;\;\;
u_{jq_jn}'=e^{m+l_j\delta} z_{jq_jn}'.
\end{equation}
Note that $u_{jsi}'$ are small at $i\leq n-1$, as $m$ is assumed to be large, and $l_j\delta$ is bounded.
The values
of $z_{jsi}'$ with $i\neq n-1$ will be taken sufficiently small as well, and we will keep
\begin{equation}\label{unz}
u_{jsn}'=u_{js} \;\; \mbox{ and } \;\; z_{j,s+1,n-1}'=-e^{-\lambda\delta} u_{js}^{-\lambda},
\end{equation}
in order to ensure the closeness of $P_{js}'$ to $P_{js}$ and $M_{js}'$ to $M_{js}$.

The first of the small perturbations which we add to the map $Y_\delta$ is
localized in a small neighborhood of the points
$Y_\delta^{-1}(P_{j,s+1})$ (so outside these small neighborhoods $Y_\delta$ remains unchanged).
We take these localized perturbations such that in a sufficiently
small neighborhood of $M_{js}=Y_{\delta}^{-N}(Y_\delta^{-1}P_{j,s+1})$ the map $\tilde Y_{\delta}^{N+1}$
(where $\tilde Y_{\delta}$ denotes the perturbed map) is given by (\ref{ydmap})
with the following correction term
\begin{equation}\label{pcr}
\begin{array}{l}\displaystyle
z_{j,s+1,n}' -
(-1)^{n+1} e^{\delta} x_n^{2-(n-2)\lambda}
\left|\frac{x_{n-1}-b_j}{z_{js,n-1}'}\right|^{\frac{1}{\lambda}-(n-1)} e^{-\lambda m} z_{js1}'+\\
\\ \displaystyle
\qquad\qquad+\;\sum_{\stackrel{\nu_2\geq 0,\dots,\nu_n\geq 0}{\nu_2+\dots+\nu_n \leq d}}
\varepsilon_{js\nu} \prod_{2\leq p\leq n} (\hat x_p-u_{jsp}')^{\nu_p}\end{array}
\end{equation}
added into the equation for $\bar x_n$, where $\varepsilon_{js\nu}$ are small coefficients
to be determined later (see (\ref{epj})).
The first term in (\ref{pcr}) is small as well (see (\ref{uzr}),(\ref{uzrl})); in the second term
the values of $x_n$ and $z_{j,s,n-1}'$ are bounded away from zero, and the exponent
$\left(\frac{1}{\lambda}-(n-1)\right)$ is larger than $r$ (see (\ref{gamd})),
hence the second term is also small with the derivatives up to the order $r$ at least.
The first two terms in (\ref{pcr}) ensure, in particular,
that at $\varepsilon=0$ the coordinate $x_n$ of $\tilde Y_\delta^{N+1} M_{js}'$
coincides with that of $P_{j,s+1}'$ (see (\ref{uzr}),(\ref{unz}),(\ref{gam})).
We want $P_{j,s+1}'=\tilde Y_\delta^{N+1} M_{js}'$  at $\varepsilon=0$, so we put
\begin{equation}\label{uzrr}
\begin{array}{l}\displaystyle
z_{j,s+1,1}'=e^{-(1-(n-2)\lambda)\delta-\lambda m} u_{js}^{(n-1)\lambda-1} z_{js2}',\\
\\ \displaystyle z_{j,s+1,i}'=e^{-\lambda (m+\delta)} z_{js,i+1}' \;\; (2\leq i\leq  n-2),
\end{array}
\end{equation}
(see (\ref{ydmap}),(\ref{uzr}),(\ref{unz})). At $s=1$ we assume
\begin{equation}\label{z1}
z_{j1i}'=0 \; \mbox{ at } \; i\leq n-2, \;\;\; z_{j1,n-1}'=e^{-\lambda\delta/2}.
\end{equation}
Now, the values of $z_{jsi}', u_{jsi}'$ are defined by (\ref{uzr}),(\ref{uzrl}),(\ref{unz}),(\ref{uzrr})
for all $j,s,i$. As one can see, $z_{jsi}'$ at $i\neq n-1$ and $u_{jsi}'$ at $i\neq n$ tend to zero
as $m\rightarrow +\infty$, i.e. $P_{js}'\rightarrow P_{js}$ and $M_{js}'\rightarrow M_{js}$ indeed.

At all small $\varepsilon$ the map
$\tilde T_{j,1+\delta}\circ L_{jm}\equiv \tilde Y_{\delta}^{mN+N+1}$ takes a small
neighborhood of $P_{js}'$ into a small neighborhood of $P_{j,s+1}'$.
We choose some $\eta(m)$ that tends to zero as $m\rightarrow+\infty$ and some, independent of $m$,
coefficients $C_{jsi}>0$, and
introduce rescaled coordinates $v_1,\dots,v_n$ near $P_{js}'$ by the rule
\begin{equation}\label{scal}
\begin{array}{l}\displaystyle
x_1 |b_j-x_{n-1}|^{n-1-1/\lambda}=
z_{js1}'|z_{js,n-1}'|^{n-1-1/\lambda}+ C_{js1} \; \eta \; e^{-\lambda m (n-2)}v_1\\ \\
\displaystyle \hat x_i=z_{jsi}'+ C_{jsi} \; \eta \; e^{-\lambda m (n-i-1)} v_i \;\;\; (2\leq i\leq n-1),
\qquad x_n=z_{jsn}'+ C_{jsn} \; \eta \; e^{-m} v_n\end{array}
\end{equation}
(recall that $|b_j-x_{n-1}|$ is close to $1$ near $P_{js}$,
hence (\ref{scal}) is a smooth coordinate transformation).
Since $\eta$ tends to zero as $m\rightarrow+\infty$, any bounded region of values
of $v$ corresponds to a small neighborhood of $P_{js}'$.

After the rescaling, the map $\tilde T_{j,1+\delta}\circ L_{jm}\equiv \tilde Y_{\delta}^{mN+N+1}$
from a small neighborhood of $P_{js}'$ into a small neighborhood of $P_{j,s+1}'$ takes the following form
(see
(\ref{ydmap}),(\ref{ltmap}),(\ref{gam}),(\ref{uzr}),(\ref{unz}),(\ref{uzrr}),(\ref{scal})):
$$\begin{array}{l} \displaystyle
C_{j,s+1,i} \bar v_i = e^{-\lambda\delta} C_{j,s,i+1} v_{i+1} \;\;\; (i\leq  n-2),\\ \\ \displaystyle
C_{j,s+1,n-1} \bar v_{n-1}=
e^{-\lambda\delta}\;
(u_{js}^{-\lambda}-(u_{js}+ C_{jsn} \;\eta\; v_n)^{-\lambda})/\eta,\\ \\
\displaystyle
C_{j,s+1,n} \bar v_n =
(-1)^{n+1} \phi_{js} C_{js1} v_1
+\;\sum_{\stackrel{\nu_2\geq 0,\dots,\nu_n\geq 0}{\nu_2+\dots+\nu_n \leq d}}
\varepsilon_{js\nu} E_{js\nu} \prod_{2\leq p\leq n} v_p^{\nu_p}
\end{array}
$$
where we denote
$$\begin{array}{l}\displaystyle
\phi_{js}=e^{\delta} (u_{js}+\eta C_{jsn} v_n)^{2-(n-2)\lambda}
|z_{js,n-1}'+ \eta C_{j,s,n-1} v_{n-1}|^{\frac{1}{\lambda}-(n-1)},\\ \\ \displaystyle
E_{js\nu}=e^{m (1-\lambda \sum_{2\leq p\leq n-1} (n-p)\nu_p)}\;
\eta^{(-1+\sum_{2\leq p\leq n} \nu_p)} \; \prod_{2\leq p\leq n} C_{jsp}^{\nu_p}.
\end{array}
$$
Note that $\sum_{2\leq p\leq n-1} (n-p)\nu_p\leq (n-2)d$, hence
$1-\lambda \sum_{2\leq p\leq n-1} (n-p)\nu_p>0$ (see (\ref{gamd})). Therefore, all
the coefficients $E_{js\nu}$ tend to infinity as $m\rightarrow+\infty$ (provided
$\eta$ tends to zero sufficiently slowly).

As we see, by putting
\begin{equation}\label{cj}
\begin{array}{l}\displaystyle
C_{j,s+1,i}=e^{-\lambda\delta} C_{j,s,i+1} \;\;\; (i\leq  n-2), \qquad
C_{j,s+1,n-1}=\lambda e^{-\lambda\delta} u_{js}^{-\lambda-1} C_{jsn},\\ \\
\displaystyle
C_{j,s+1,n}=e^{\delta} u_{js}^{2-(n-2)\lambda}
|z_{js,n-1}'|^{\frac{1}{\lambda}-(n-1)} C_{js1},
\end{array}
\end{equation}
and
\begin{equation}\label{epj}
\varepsilon_{js\nu}=h_{jsi} \frac{C_{j,s+1,n}}{E_{js\nu}},
\end{equation}
the map $\tilde T_{j,1+\delta}\circ L_{jm}$ near $P_{js}'$ takes the form
\begin{equation}\label{scsdl}
\begin{array}{l} \displaystyle
\bar v_i = v_{i+1} \;\;\; (i\leq  n-2),\qquad
\bar v_{n-1}= v_n + O(\eta),\\ \\
\displaystyle
\bar v_n = (-1)^{n+1} v_1 + O(\eta)
+\;\sum_{\stackrel{\nu_2\geq 0,\dots,\nu_n\geq 0}{\nu_2+\dots+\nu_n \leq d}}
h_{js\nu} \prod_{2\leq p\leq n} v_p^{\nu_p},
\end{array}
\end{equation}
i.e. it can be made as close as we want to the map $\tilde H_{js}$, provided $m$ is taken large enough
(recall that $\eta\rightarrow 0$ as $m\rightarrow+\infty$). We take $\eta$ tending to zero sufficiently
slowly, so, as we mentioned, $E_{ij\nu}\rightarrow \infty$ as $m\rightarrow+\infty$, which implies
that all $\varepsilon_{js\nu}\rightarrow 0$ (see (\ref{epj})), i.e. our perturbation to $Y_\delta$
is arbitrarily small indeed.

It follows that in the rescaled coordinates the map
$\left(\tilde T_{j,1+\delta}\circ L_{jm}\right)^{q_j}\equiv \tilde Y_\delta^{q_j(mN+N+1)}$
from a small neighborhood of $P_{j1}'$ into a small neighborhood of $P_{jq_j}'$
can be made as close as we want to the map $\tilde H_{jq_j}\circ \dots \circ \tilde H_{j1}$,
provided $m$ is large enough (the rescaled coordinates near $P_{j1}'$ and $P_{jq_j}'$ are
given by formulas (\ref{scal}), where the coefficients $C_{j1i}>0$ are taken arbitrary,
and the coefficients $C_{jq_ji}$ are recovered from the recursive formula (\ref{cj}); since
$m$ does not enter (\ref{cj}), it follows that $C_{jq_ji}$ stay bounded away from zero
and infinity as $m\rightarrow+\infty$).

Now, from (\ref{ltmap}) we obtain that the same holds true for the map
$L_{j,m+l_j\delta}\circ\left(\tilde T_{j,1+\delta}\circ L_{jm}\right)^{q_j}\equiv
\tilde Y_\delta^{l_j+mN+q_j(mN+N+1)}$
from a small neighborhood of $P_{j1}'$ into a small neighborhood of $M_{jq_j}'$,
where the rescaled coordinates $(v_1,\dots,v_n)$ are introduced as follows:
\begin{equation}\label{scall}
\begin{array}{l}\displaystyle
x_1 |x_{n-1}-b_j|^{n-1-1/\lambda}=u_{jq_j1}'|u_{jq_j,n-1}'|^{n-1-1/\lambda}+
C_{jq_j1}\; \eta \; e^{-\lambda (n-1) m-\lambda l_j\delta} v_1,\\ \\
\displaystyle \hat x_i=u_{jq_ji}'+ C_{jq_ji} \; \eta \; e^{-\lambda m (n-i)-\lambda l_j\delta} v_i
\;\; (2\leq i\leq n-1),\qquad
\displaystyle x_n=u_{jq_j} + C_{jq_jn} \; \eta \; e^{l_j\delta} v_n,
\end{array}
\end{equation}
with the same constants $C_{jq_ji}$ as above.

Recall that, by construction, the point $Y_\delta^{N+1}M_{jq_j}'$ lies
in $U_{3-j,+}$ in the region $a_{3-j,+}+1-\delta<x_{n-1}<a_{3-j,+}+1$. We add to
the map $Y_\delta$ an additional perturbation, localized near the point $Y_\delta^{N}M_{jq_j}'$, such that
the corresponding map $\tilde G_{j,1+\delta}\equiv \tilde Y_\delta^{N+1}$ will have the following form
near $M_{jq_j}'$:
\begin{equation}\label{gtmapt}
\begin{array}{l} \displaystyle
\bar x_1=|x_n|^{1-(n-2)\lambda}
(x_1-e^{-\lambda (m+l_j\delta)} z_{jq_j1}'\left|(x_{n-1}-b_j)/z_{jq_j,n-1}'\right|^{\frac{1}{\lambda}-(n-1)}),
\\ \\ \displaystyle
\bar x_i=|x_n|^{\lambda} (x_i-u_{jq_ji}') \;\; (2\leq i\leq  n-2),\\ \\ \displaystyle
\bar x_{n-1}=a_{3-j,+}+1-\delta-\ln |x_n|,\qquad
\bar x_n=-|x_{n}|^{\lambda} (x_{n-1}-b_j-u_{jq_j,n-1}').
\end{array}
\end{equation}
Note that $u_{jq_ji}'$ at $i\leq n-1$ tend to zero as $m\rightarrow+\infty$ (see (\ref{uzrl})),
while the values of $x_n$ near $M_{jq_j}'$ and $z_{j,q_j,n-1}'$ are bounded away
from zero; the exponent $\frac{1}{\lambda}-(n-1)$ in the first line is larger than $r$ (see (\ref{gamd})).
Thus, for sufficiently large $m$, map (\ref{gtmapt}) is indeed a small perturbation of the map
$G_{j,1+\delta}$ given by (\ref{gtmap}).

Denote $P_{3-j,0}'=\tilde G_{j,1+\delta} M_{jq_j}'$. By (\ref{gtmapt}), this is the point with the
coordinates $x_i=0$ at $i\neq n-1$, and $x_{n-1}=a_{3-j,+}+1+(\kappa_j-1)\delta$ (we assume that
the coordinate $x_n$ of $M_{jq_j}'$ is $u_{jq_jn}'=u_{jq_j}=-e^{-\kappa_j \delta}$
where $\kappa_j\in(0,1]$ is defined by (\ref{ljdf})). Introduce rescaled coordinates near $P_{3-j,0}'$ by the rule
\begin{equation}\label{resp3}
\begin{array}{l}\displaystyle
\displaystyle x_i=e^{-(l_j+\kappa_j)\delta\lambda_i}C_{jq_ji} \; \eta \;
e^{-\lambda m (n-i)} v_i \;\; (i\leq n-2),\\
\displaystyle x_{n-1}=a_{3-j,+}+1+(1-\kappa_j)\delta+e^{(l_j+\kappa_j)\delta} C_{jq_jn} \; \eta \; v_{n-1},\\
\displaystyle x_n=e^{-(l_j+\kappa_j)\delta\lambda} C_{j,q_j,n-1} \; \eta \; e^{-\lambda m} v_n
\end{array}
\end{equation}
(with the same constants $C_{jq_ji}$ as above).
In coordinates (\ref{scall}),(\ref{resp3}), map (\ref{gtmapt}) takes the form
$(v_1,\dots,v_{n-1},v_n)\mapsto(v_1,\dots,v_n,-v_{n-1})+O(\eta)$, i.e. it becomes arbitrarily
close to the map $\Phi_0$ (see (\ref{phi0})) as $m\rightarrow+\infty$. Thus, in the
rescaled coordinates, the map
$\tilde G_{j,1+\delta}\circ L_{j,m+l_j\delta}\circ\left(\tilde T_{j,1+\delta}\circ L_{jm}\right)^{q_j}\equiv
\tilde Y_\delta^{l_j+(1+q_j)(mN+N+1)}$
from a small neighborhood of $P_{j1}'$ into a small neighborhood of $P_{3-j,0}'$ ,
can be made as close as we want to the map $\Phi_0\circ\tilde H_{jq_j}\circ \dots \circ \tilde H_{j1}$
as $m$ grows.

Analogously, we take the point $M_{j0}':\{x_{n-1}=a_{j-}-1+\delta/2, x_i=0 \; (i\neq n-1)\}\in U_{j-}$,
and perturb the map $Y_\delta$ near $Y_\delta^N P_{j0}'$ in such a way that the map
$\tilde Q_{j,1+\delta}\equiv \tilde Y_\delta^{N+1}$ near $M_{j0}'$ will be given by
\begin{equation}\label{qtmapm}
\begin{array}{l} \displaystyle
\bar x_i=e^{-\lambda_i(\delta-x_{n-1}-1+a_{j-})}x_i \qquad (i\leq  n-2),\\ \\ \displaystyle
\bar x_{n-1}-b_j=e^{-\lambda(\delta-x_{n-1}-1+a_{j-})},\;\;
\bar x_n=e^{\delta-x_{n-1}-1+a_{j-}} x_n + e^{-m}u_{j1}.
\end{array}
\end{equation}
It is a small perturbation of the map $Q_{j,1+\delta}$ from (\ref{qtmap}), and it takes
$M_{j0}'$ to $P_{j1}'$ (see (\ref{uzr}),(\ref{unz})). When we introduce rescaled variables
near $M_{j0}'$ by the rule
\begin{equation}\label{resp0}
\begin{array}{l}
\displaystyle x_i=e^{\delta\lambda_i/2}C_{j1i} \; \eta \; e^{-\lambda m (n-i-1)} v_i \;\; (i\leq n-2),\\
\displaystyle x_{n-1}=a_{j-}-1+\delta/2+
\frac{1}{\lambda}e^{\lambda\delta/2} C_{j,1,n-1} \; \eta\; v_{n-1},\\
\displaystyle x_n=e^{-\delta/2} C_{j1n} \; \eta \; e^{-m} v_n,
\end{array}
\end{equation}
map (\ref{qtmapm}) will take the form $\bar v=v+O(\eta)$, i.e. it is close to the identity map.
Thus, in the rescaled coordinates given by (\ref{resp0}),(\ref{resp3}), the map
$\tilde G_{j,1+\delta}\circ L_{j,m+l_j\delta}\circ
\left(\tilde T_{j,1+\delta}\circ L_{jm}\right)^{q_j}\circ\tilde Q_{j,1+\delta} \equiv
\tilde Y_\delta^{l_j+(1+q_j)(mN+N+1)+N+1}$
from a small neighborhood of $M_{j0}'$ into a small neighborhood of $P_{3-j,0}'$,
is as close as we want to the map $\Phi_0\circ\tilde H_{jq_j}\circ \dots \circ \tilde H_{j1}$
at $m$ large enough, i.e. it is a close approximation of the map $\Phi_j$.

Let us now determine the form of the map $S_{jt}:v\mapsto\bar v$ from a small neighborhood of $P_{j0}'$
into
a small neighborhood of $M_{j0}'$ in the rescaled coordinates (\ref{resp0}),(\ref{resp3}).
By (\ref{stmap}), for an integer $k>0$, the map $S_{j,k\delta}\equiv Y_{\delta}^{k}$ takes the point
$P_{j0}'$ into $M_{j0}'$ if
\begin{equation}\label{mdf}
\beta(\mu_j)= (k+\kappa_{3-j}-1/2)\delta-5
\end{equation}
(see (\ref{resp0}),(\ref{resp3})).
Since $\beta\rightarrow+\infty$ as $\mu\rightarrow+0$ (see (\ref{albet})),
for every
sufficiently large $k$ equation (\ref{mdf}) has a solution $\mu_j(k)$, and $\mu_j(k)\rightarrow +0$
as $k\rightarrow+\infty$. It follows that $\alpha(\mu_j(k))\rightarrow+\infty$. Thus, for any sufficiently
large $m$ we can find $\gamma_{n\pm}\in (0,1]$ and $k$ such that
\begin{equation}\label{mdfa}
\begin{array}{l}\displaystyle
e^{-\gamma_{_{n+}}\alpha(\mu_j(k))} =
e^{-(l_{3-j}+\kappa_{3-j})\delta\lambda} C_{3-j,q_{3-j},n-1} \; \eta \; e^{-\lambda m},\\ \\
\displaystyle
e^{-\gamma_{_{n-}}\alpha(\mu_j(k))} = e^{-\delta/2} C_{j1n} \; \eta \; e^{-m}.
\end{array}
\end{equation}
This guarantees that
$\bar v_n=\psi_j(v_n)$ (see (\ref{resp0}),(\ref{resp3}),(\ref{stmap})).

We also obtain $\bar v_i=v_i$ at $i\leq n-2$ by choosing $\gamma_{i\pm}\in(0,1]$ such that
\begin{equation}\label{gdf}
\begin{array}{l}\displaystyle
e^{-\gamma_{_{i+}}\alpha(\mu_j(k))} =
e^{\delta\lambda_i/2}C_{j1i} \; \eta \; e^{-\lambda m},\\ \\ \displaystyle
e^{-\gamma_{_{i-}}\alpha(\mu_j(k))} = e^{-(l_j+\kappa_j)\delta\lambda_i}C_{3-j,q_{3-j},i}\eta.
\end{array}
\end{equation}

Finally, we fix the choice of the integer $l_j$ and $\kappa_j\in(0,1]$ as follows:
\begin{equation}\label{ljdf}
e^{(l_j+\kappa_j)\delta} =\frac{1}{\lambda}e^{\lambda\delta/2} C_{3-j,1,n-1}/C_{jq_jn}.
\end{equation}
This (along with (\ref{mdf})) gives us $\bar v_{n-1}=v_{n-1}$ for the map $S_{j,k\delta}$
in the coordinates (\ref{resp0}),(\ref{resp3}). As wee see,
the map $S_{j,k\delta}$ in the rescaled coordinates coincides with the map $\Psi_j$
for $v$ from some open neighborhood of $D$ (if $j=1$) or of $\Phi_1\circ\Psi_1(D)$ (if $j=2$).

Thus, we see that the map
$$\begin{array}{l} \displaystyle
\tilde G_{2,1+\delta}\circ L_{2,m+l_2\delta}\circ
\left(\tilde T_{2,1+\delta}\circ L_{2m}\right)^{q_2}\circ\tilde Q_{2,1+\delta}\circ S_{2,k\delta}\circ\\
\displaystyle
\qquad\qquad \circ\; \tilde G_{1,1+\delta}\circ L_{1,m+l_1\delta}\circ
\left(\tilde T_{1,1+\delta}\circ L_{1m}\right)^{q_1}\circ\tilde Q_{1,1+\delta}\circ S_{1,k\delta}\equiv\\
\displaystyle \equiv \;
\tilde Y_\delta^{2k+l_1+l_2+(2+q_1+q_2)(mN+N+1)+2(N+1)}
\end{array}$$
is a close approximation to the map $F=\Phi_2\circ\Psi_2\circ\Phi_1\circ\Psi_1$,
provided $\tilde H_{jq_j}\circ\dots \circ \tilde H_{j1}$ are sufficiently close approximations to
$\tilde \Phi_j$ ($j=1,2$) and $m$ is large enough. This completes the proof of the theorem.

\section{Birth of periodic spots from a heteroclinic cycle}
\label{henhom}\setcounter{equation}{0}
\renewcommand{\theequation}{\thesection.\arabic{equation}}
In this Section we finish the proof of Theorem \ref{thmn}.
Let a $C^\rho$-diffeomorphism $f$
of a smooth two-dimensional manifold have a pair of saddle periodic points
$P$ and $Q$ of periods $p$ and, respectively, $q$. Denote as $T_{01}$ the map $f^p$ restricted onto a small
neighborhood of $P$, and denote as $T_{02}$ the map $f^q$ restricted onto a small neighborhood of $Q$.
By definition, $T_{01}P=P$ and $T_{02}Q=Q$.
One can introduce coordinates $(x_1,y_1)$ in the neighborhood of $P$ and $(x_2,y_2)$ in the neighborhood
of $Q$ such that the maps $T_{0j}$ will have the form
$$\bar x_j=\lambda_j x_j +\dots, \qquad \bar y_j=\gamma_j y_j +\dots,$$
where $|\lambda_j|<1$, $|\gamma_j|>1$; the dots stand for nonlinearities. The numbers $\lambda_1,\gamma_1$
and $\lambda_2,\gamma_2$ are the multipliers of the periodic points $P$ and $Q$, respectively.
We assume
\begin{equation}\label{sigmas}
J_1:=|\lambda_1\gamma_1|<1 \;\;\;\mbox{ and }\;\;\;J_2:=|\lambda_2\gamma_2|>1.
\end{equation}

\begin{figure}
\centerline{\includegraphics[scale=0.5]{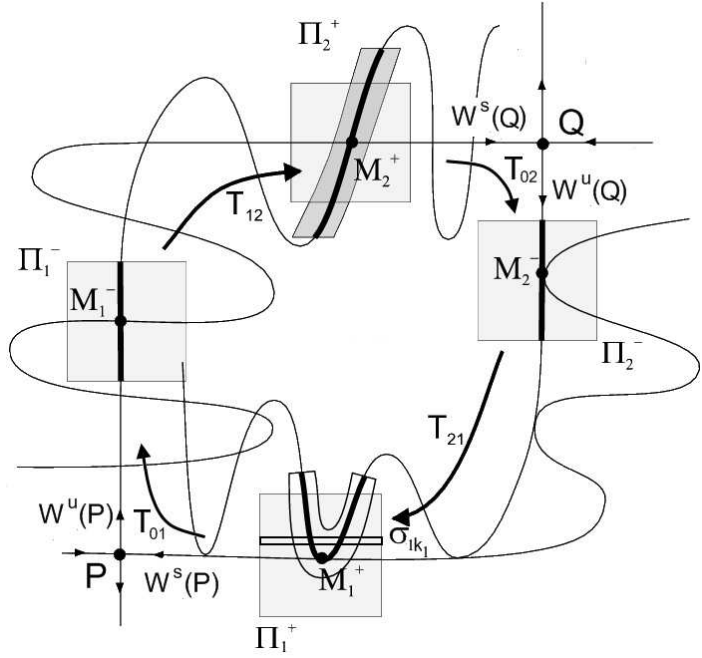}}
\caption{A non-transverse heteroclinic cycle}
\end{figure}

Every saddle periodic point lies in the intersection of two $C^\rho$-smooth invariant manifolds: the points
in the stable invariant manifold $W^s$ tend to the periodic orbit at forward iterations of the map, and the
points in the unstable invariant manifold $W^u$ tend to the periodic orbit at backward iterations. With our
choice of the coordinates, the unstable manifold is tangent at the periodic point to the $y$-axis, and
the stable manifold is tangent to the $x$-axis. One can locally straighten $W^u$ and $W^s$, i.e.
the coordinates $(x_j,y_j)$ can be chosen in such a way that $W^s_{loc}(P)=\{y_1=0\}$,
$W^u_{loc}(P)=\{x_1=0\}$, $W^s_{loc}(Q)=\{y_2=0\}$, $W^u_{loc}(Q)=\{x_2=0\}$.

Assume that $W^u(P)$ has an orbit of a transverse intersection with $W^s(Q)$. This means that
$M_2^+:=f^{k_{12}}M_1^-\in W^s_{loc}(Q)$ for some positive integer $k_{12}$ and some point
$M_1^-(0,y_1^-)\in W^u_{loc}(Q)$, and the curve $f^{k_{12}}(W^u_{loc}(P))$ intersects $W^s_{loc}(Q)$
at the point $M_2^+(x_2^+,0)$ transversely. We denote the map $f^{k_{12}}$ restricted onto a small
neighborhood of $M_1^-$ as $T_{12}$. It can be written as
\begin{equation}\label{t12}
\bar x_2 - x_2^+ = a_1 x_1 +b_1 (y_1-y_1^-)+\dots, \qquad \bar y_2=c_1 x_1 +d_1 (y_1-y_1^-)+\dots,
\end{equation}
where $d_1\neq 0$ because of the transversality of $T_{12}(W^u_{loc}(P))$ to $W^s_{loc}(Q)$. As an orbit of
intersection of $W^u(P)$ and $W^s(Q)$, the orbit $\Gamma_{PQ}$ of the point $M_1^-$ is heteroclinic:
it tends to the orbit of $P$ at backward iterations of $f$ and to the orbit of $Q$ at forward iterations.
Note that since $d_1\neq0$, we may rewrite (\ref{t12}) in the so-called cross-form:
\begin{equation}\label{t12c}
d_1(\bar x_2 - x_2^+) = D x_1 + b_1 \bar y_2 +\dots, \qquad d_1 (y_1-y_1^-)=\bar y_2 -c_1 x_1 +\dots,
\end{equation}
where $D:=a_1d_1-b_1c_1$.

Another assumption is that $f$ has a heteroclinic orbit $\Gamma_{QP}$ at the points of which $W^u(Q)$ has
a tangency of order $m$ with $W^s(P)$. This means that there exists a pair of points,
$M_2^-(0,y^-_2)\in W^u_{loc}(Q)$ and $M_1^+(x_1^+,0)\in W^s_{loc}(P)$, such that $M_1^+=f^{k_{21}}M_2^-$
for some positive integer $k_{21}$, and the curve $f^{k_{21}}(W^u_{loc}(Q))$ has a tangency of order $m$
with $W^s_{loc}(P)$ at $M_1^+$ (see Fig.3). We denote the map $f^{k_{21}}$ restricted onto
a small neighborhood of $M_2^-$ as $T_{21}$. It can be written as
\begin{equation}\label{t21}
\bar x_1 - x_1^+ = a_2 x_2 +b_2 (y_2-y_2^-)+\dots, \qquad \bar y_1=c_2 x_2 +d_2(y_2-y_2^-)^{m+1}+\dots,
\end{equation}
where $d_2\neq 0$; the dots stand for higher order terms. Obviously, to speak about the tangency
of order $m$, the smoothness of $f$ has to be sufficiently high, i.e. $\rho\geq m+1$.

We will fix the orientation in the neighborhood of $Q$ by requiring that the determinant $D$
of the derivative matrix $\partial T_{12}/\partial(x,y)$ is positive at the point $M^-_1$.
Then, we require that the
determinant of $\partial T_{21}/\partial(x,y)$ at the point $M^-_2$ is also positive, i.e. $b_2c_2<0$.
This is always the case if the manifold on which $f$ is defined is orientable and $f$
is orientation-preserving. In Remark after Lemma \ref{thmn2} we discuss the case $b_2c_2>0$ too.

In Lemma \ref{thmn2} we will also need a technical assumption
$\ln|\gamma_1|\ln|\gamma_2|\neq\ln|\lambda_1|\ln|\lambda_2|$ (one can always achieve this by an arbitrarily small perturbation of $f$ without
destroying the order $m$ tangency between $W^u(Q)$ and $W^s(P)$). In fact, we
can always assume
\begin{equation}\label{tc}
\ln|\gamma_1|\ln|\gamma_2| < \ln|\lambda_1|\ln|\lambda_2|;
\end{equation}
the case $\ln|\gamma_1|\ln|\gamma_2| > \ln|\lambda_1|\ln|\lambda_2|$ reduces to the given one by
considering the map $f^{-1}$ instead of $f$ (and interchanging $P$ with $Q$).

Let us imbed $f$ in any $(m+1)$-parameter family $f_\nu$ of maps, of class $C^\rho$ with
respect to coordinate and parameters. The corresponding map $T_{21}$ will also depend smoothly
on parameters. Since it has form (\ref{t21}) at $\nu=0$, i.e. all the derivatives of $\bar y_1$
with respect to $(y_2-y_2^-)$ vanish up to the order $m$, it follows that at non-zero $\nu$
the map $T_{21}$ can be written in the form
\begin{equation}\label{t21m}
\bar x_1 - x_1^+ = a_2 x_2 +b_2 (y_2-y_2^-)+\dots, \qquad \bar y_1=c_2 x_2 +
\sum_{s=0}^{m-1} \mu_s (y_2-y^-_2)^s + d_2(y_2-y_2^-)^{m+1}+\dots,
\end{equation}
where $\mu_s$ are smooth functions of $\nu$ (note that $x_1^+$ and $y_2^-$ also depend on $\nu$ now:
the value of $y_2^-$ is fixed by the condition $\partial^m \bar y_1/\partial y_2^m=0$ at $y_2=y_2^-$;
thus $y_2^-$ is $C^{\rho-m}$-function of $\nu$, and $\mu_s$ are also $C^{\rho-m}$; the high-order terms
that are denoted by dots in (\ref{t21m}) and (\ref{t12c}) depend now also on $\nu$,
$C^{\rho-m}$-smoothly).

Denote $\theta=|\ln J_2/\ln J_1|$.
The value of $\theta$ is a $C^{\rho-1}$-smooth function of $\nu$.
We may always put a family $f_\nu$ in general position, i.e. we may further assume
\begin{equation}\label{gennu}
\frac{\partial(\mu_0,\dots,\mu_{m-1},\theta)}{\partial(\nu_1,\dots,\nu_{m+1})}\neq0.
\end{equation}
This, in particular, means that by changing $\nu$ we may change the values of any of the
parameters $\theta$ and $\mu_s$ while keeping the other parameters constant. For example,
one may change the value of $\theta$ and keep $\mu=0$, i.e. keep the order $m$ tangency between
$W^u(Q)$ and $W^s(P)$.

The two periodic orbits (the orbit of $P$ and the orbit of $Q$) and two heteroclinic orbits,
$\Gamma_{PQ}$ and $\Gamma_{QP}$ comprise a heteroclinic cycle. Let $U$ be a small neighborhood of
the heteroclinic cycle. Varying the parameters $\nu$ can lead to the destruction of the heteroclinic
cycle and to bifurcations of other orbits in $U$. We will further focus on one instance of such
bifurcations. We call a periodic orbit of the map $f_\nu$ {\em single-round} if
it stays entirely in $U$ and
visits a small neighborhood of each of the points $M^{\pm}_{1,2}$ only once, i.e. the point of
intersection of the single-round orbit with a small neighborhood of $M_1^+$ is a fixed point
of the map $T_{21}T_{02}^{k_2}T_{12}T_{01}^{k_1}$ for some positive integers $k_1$ and $k_2$; the image
of this point by the map $T_{01}^{k_1}$ lies in a small neighborhood of $M^+_1$, the image
by $T_{12}T_{01}^{k_1}$ lies in a small neighborhood of $M^-_2$, and the image by
$T_{02}^{k_2}T_{12}T_{01}^{k_1}$ lies in a small neighborhood of $M^+_2$.

\begin{lemma}$\!\!.$\label{thmn2} There exist a sequence $\nu_l\rightarrow 0$ and sequences of integers
$k_{1l}\rightarrow+\infty$ and $k_{2l}\rightarrow+\infty$ such that at $\nu=\nu_l$ the map $f_\nu$
has a single-round periodic orbit which corresponds to a fixed point of the map
${\cal T}_l:=T_{21}T_{02}^{k_{2l}}T_{12}T_{01}^{k_{1l}}$, and the map ${\cal T}_l$ near this
point is, in some $C^{\rho}$-coordinates $(u,v)$, given by
\begin{equation}\label{flat}
(\bar u,\bar v)=\Phi(u,v)+o(|u|^m+|v|^m),
\end{equation}
where $\Phi$ denotes the time-$1$ map by the flow
\begin{equation}\label{flowat}
\dot u=v,\qquad \dot v = -\Psi(u)(1+v)
\end{equation}
near $(u,v)=(0,0)$; here $\Psi$ is a polynomial such that $\Psi(0)=0$ and $\Psi'(0)\geq0$.
\end{lemma}
\noindent{Remark.} System (\ref{flowat}) has an integral
$H(u,v)=\int\Psi(u)du + v -\ln(1+v)=\frac{v^2}{2}+\Psi'(0)\frac{u^2}{2}+\dots$. Thus, when $\Psi'(0)>0$,
the equilibrium at zero is a center: every orbit of (\ref{flowat}) is in this case a closed curve
surrounding the origin. The corresponding time-$1$ map $\Phi$ is, therefore, conservative (it preserves
the area form $\frac{1}{1+v}\;du\wedge dv$) and has an elliptic point at the origin. Note that map
(\ref{flat}) can, by an arbitrarily $C^m$-small perturbation, be made equal to
$(\bar u,\bar v)=\Phi(u,v)$  identically
in a sufficiently small neighborhood of zero. Hence, Lemma \ref{thmn2}
implies, that by an arbitrarily $C^m$-small perturbation of the given map $f$
a periodic point can be born in a small neighborhood of the heteroclinic cycle such that the first-return
map near this point will be area-preserving and the point will be elliptic. We recall that
one of the conditions of the lemma is that the heteroclinic cycle is ``orientable'' in the sense that
the determinant of $\partial T_{21}/\partial(x,y)$ at the point $M^-_2$ is positive. However, it is easy
to show (see Lemma 8 in \cite{GTS6}) that a non-orientable heteroclinic cycle can always be
perturbed in such a way that a new (double-round) orbit of heteroclinic tangency of order $(m-1)$
between $W^u(Q)$ and $W^s(P)$ is born, and the corresponding cycle is now orientable; applying the lemma
to the newly born cycle, we find that an elliptic periodic point can be born by a perturbation which
is arbitrarily small in the $C^{m-1}$-metric. In any case, we have that for any fixed $r$, given
a $C^\infty$-map $f$ which satisfies (\ref{sigmas}) and which has a heteroclinic cycle with a
sufficiently high order of tangency between $W^u(Q)$ and $W^s(P)$, an arbitrarily $C^r$-small
perturbation leads to the birth of an elliptic periodic point: a point whose both multipliers lie
on the unit circle and, importantly, the corresponding first-return map is area-preserving.
By \cite{MR}, an arbitrarily $C^r$-small perturbation of the map in a neighborhood of such point
creates an orbit of homoclinic tangency; the perturbation does not destroy the area-preserving
property. By Theorem 5 of \cite{GTS6}, an arbitrarily $C^r$-small perturbation of any area-preserving
map with a homoclinic tangency leads to the birth of a periodic spot. In other words, once
Lemma \ref{thmn2} provides us with a periodic orbit whose normal form is area-preserving up to
a sufficiently high order, the birth of periodic spots from the heteroclinic cycles under
consideration follows from the results of \cite{GTS6} on perturbations of area-preserving maps.\\

{\em Proof of Lemma \ref{thmn2}.} At $m=1$ an equivalent
statement can be found in \cite{GTS97,St}, so we further
focus on the case $m\geq2$. Since the periodic points $P$ and $Q$ are saddle, it follows that
given any small $x_{j0}$ and $y_{jk}$ (where $j=1,2$) and any $k\geq 0$ there exist uniquely defined small
$x_{jk}$ and $y_{j0}$ such that $(x_{jk},y_{jk}) = T_{0j}^k(x_{j0},y_{j0})$ and all the points in the
orbit $\{(x_{j0)},y_{j0}),T_{0j}(x_{j0},y_{j0}),\dots,T_{0j}^k(x_{j0},y_{j0})\}$ lie in a small
neighborhood of $P$ (at $j=1$) or $Q$ (at $j=2$), see \cite{GTS97,book}. Denote
\begin{equation}\label{xieta}
x_{jk}=\lambda_j^k x_{j0} +\xi_{jk}(x_{j0},y_{jk},\nu), \quad
y_{j0}=\gamma_j^{-k}y_{jk} + \eta_{jk}(x_{j0},y_{jk},\nu).
\end{equation}
By \cite{book}, one can introduce $C^\rho$-coordinates $(x_j,y_j)$
near the saddle periodic points in such a way that
\begin{equation}\label{xee}
\xi_{jk}=o(\lambda_j^k), \quad \eta_{jk} = o(\gamma_j^{-k}),
\end{equation}
i.e. the map $T_{0j}^k$ written in the ``cross-form'' (\ref{xieta}) is linear in the main order.
By \cite{GTS7}, the same $o(\lambda_j^k)$ and, resp., $o(\gamma_j^{-k})$ estimates hold for the
derivatives of the functions $\xi_{jk}$ and $\eta_{jk}$ up to the order $(\rho-1)$ with
respect to $(x_{j0},y_{jk})$ and up to the order $(\rho-2)$ with respect to
parameters; while for the higher order derivatives up to the order $\rho$ for
which the number of differentiations with respect to the parameters does not exceed $(\rho-2)$
we have that they uniformly tend to zero as $k\rightarrow\infty$.

Let $\delta>0$ be sufficiently small. Consider rectangular neighborhoods
$\Pi^+_j:\{|x_j-x^+_j|<\delta,|y_j|<\delta\}$ and
$\Pi^-_j:\{|y_j-y^-_j|<\delta,|x_j|<\delta\}$ of $M^+_j$ and $M^-_j$ ($j=1,2$), and take any
sufficiently large integers $k_1$ and $k_2$.
By (\ref{xieta}),(\ref{xee}), the set $\sigma_{jk_j}=\Pi^+_j\cap T^{-k_j}_{0j}\Pi^-_j$ is non-empty:
it is a strip of the form $\{|\gamma_j^{k_j} y_j - y^-_j| <\delta+ o(1)_{k\rightarrow+\infty}\}$.
We introduce a new coordinate $y_j$ on $\sigma_{jk_j}$ such that
\begin{equation}\label{zam}
y_{j,old}=\gamma_j^{-k_j}y_{j,new} + \gamma_j^{-k_j}\eta_{jk_j}(x_j,y_{j,new})
\end{equation}
(i.e. $y_{j,new}$ equals to $y_{jk_j}$ from (\ref{xieta})).

By construction, the map $T_{12}T_{01}^{k_1}$ is defined
on $\sigma_{1k_1}$; in the new coordinates, when this map takes a point $(x_1,y_1)\in \sigma_{1k_1}$ into
a point $(\bar x_2,\bar y_2)\in \sigma_{2k_2}$ it can be written in the following form (see
(\ref{t12c}),(\ref{zam})):
\begin{equation}\label{tldtk}
\begin{array}{l}\displaystyle
d_1(\bar x_2- x_2^+) = D\lambda_1^{k_1} x_1 + \phi_1(x_1) + b_1\gamma_2^{-k_2}\bar y_2 +
o(\gamma_2^{-k_2}),\\ \displaystyle
d_1(y_1-y^-_1)=- c_1\lambda_1^{k_1} x_1 - \phi_2(x_1) + \gamma_2^{-k_2}\bar y_2 + o(\gamma_2^{-k_2}),
\end{array}
\end{equation}
where $\phi_{1,2}(x_1)=o(\lambda_1^{k_1})$.
Analogously, the map $T_{21}T_{02}^{k_2}$ that takes
a point $(x_2,y_2)\in \sigma_{2k_2}$ into a point $(\bar x_1,\bar y_1)\in \sigma_{1k_1}$ is given by
\begin{equation}\label{tldtk2}
\begin{array}{l}\displaystyle
\bar x_1- x_1^+ = a_2\lambda_2^{k_2} x_2 + b_2(y_2-y^-_2) + o(y_2-y^-_2)+ o(\lambda_2^{k_2}),\\
\displaystyle
\gamma_1^{-k_1}\bar y_1=
c_2\lambda_2^{k_2} x_2 + \sum_{s=0}^{m-1} \mu_s (y_2-y_2^-)^s + d_2(y_2-y_2^-)^{m+1}+o((y_2-y_2^-)^{m+1})
+o(\lambda_2^{k_2}) + o(\gamma_1^{-k_1})
\end{array}
\end{equation}
(see (\ref{t21m})). We will now change the variable $y_1$ to
$y_1+d_1^{-1}(c_1\lambda_1^{k_1} x_1 + \phi_2(x_1))$. Then the second lines in
(\ref{tldtk}) and (\ref{tldtk2}) will change to
$$d_1(y_1-y^-_1)=\gamma_2^{-k_2}\bar y_2 + o(\gamma_2^{-k_2}),$$
and, respectively,
$$\gamma_1^{-k_1}\bar y_1=
c_2\lambda_2^{k_2} x_2 + \sum_{s=0}^{m} \hat\mu_s (y_2-y_2^-)^s + d_2(y_2-y_2^-)^{m+1}+o((y_2-y_2^-)^{m+1})
+o(\lambda_2^{k_2}) + o(\gamma_1^{-k_1}),$$
where
$\hat \mu_s=\mu_s+O(\lambda_1^{k_1})$ ($s=0,\dots,m-1$) and $\hat\mu_m=O(\lambda_1^{k_1})$.

It is easy to see, that since $d_1\neq0$, $d\neq0$, one can find constants
$C_j=O(|\lambda_{3-j}|^{k_{3-j}}+|\gamma_j|^{-k_j})$ and
$K_j=O(|\gamma_{3-j}|^{-k_{3-j}}+|\lambda_j|^{k_j})$ ($j=1,2$) such that
after the following shift of the origin:
\begin{equation}\label{shifto}
x_{j,new} = x_j - x_j^+ +C_j, \qquad y_{j,new} = y_j - y_j^- +K_j,
\end{equation}
we will have $(\bar x_{2,new},y_{1,new})=0$ at $(x_{1,new},\bar y_{2,new})=0$, and
$\bar x_{1,new}=0$, $\partial^m \bar y_{1,new}/\partial y_{1,new}^m=0$ at $(x_{2,new},y_{2,new})=0$.
Thus, after this transformation, the maps $T_{12}T_{01}^{k_1}$ and $T_{21}T_{02}^{k_2}$ will be
written as
\begin{equation}\label{tlm}
\begin{array}{l}\displaystyle
d_1 \bar x_2 = D\lambda_1^{k_1} x_1 + O(\gamma_2^{-k_2}\bar y_2) +
o((|\lambda_1|^{k_1}+|\gamma_2|^{-k_2})x_1),\\ \displaystyle
d_1 y_1= \gamma_2^{-k_2}\bar y_2 + o(\gamma_2^{-k_2}\bar y_2)+o(|\gamma_2|^{-k_2}x_1),
\end{array}
\end{equation}
and, respectively,
\begin{equation}\label{tlm2}
\begin{array}{l}\displaystyle
\bar x_1 = b_2 y_2 + o(y_2)+ O(\lambda_2^{k_2}x_2),\\ \displaystyle
\gamma_1^{-k_1}\bar y_1=
c_2\lambda_2^{k_2} x_2 + \sum_{s=0}^{m-1} \tilde\mu_s y_2^s + d_2 y_2^{m+1}+o(y_2^{m+1})
+o(\lambda_2^{k_2}x_2),
\end{array}
\end{equation}
where the modified parameters $\tilde \mu_s$ are such that
$\tilde\mu_s=\mu_s +o(1)_{k_{1,2}\rightarrow+\infty}$.

By virtue of (\ref{gennu}), by an arbitrarily small change of $\nu$ we can make
$\theta:=\left|\ln|\lambda_2\gamma_2|/\ln|\lambda_1\gamma_1|\right|$
rational without changing the values of $\mu$. Therefore,
we may from the very beginning assume that at $\nu=0$ the system has a heteroclinic cycle
with a tangency of the order $m$, and the value $\theta_0:=\theta(0)$ is rational. Further we
always assume that
\begin{equation}\label{kb}
k_1=\theta_0 k_2
\end{equation}
in (\ref{tlm}),(\ref{tlm2}). We will also assume that both $k_1$ and $k_2$ are even, so
$\lambda_j^{k_j}$ and $\gamma_j^{k_j}$ are positive.
It follows from (\ref{kb}),(\ref{tc}),(\ref{sigmas}) that
\begin{equation}\label{lgin}
\gamma_2^{-k_2}\ll\lambda_1^{k_1}\ll\gamma_1^{-k_1}\ll\lambda_2^{k_2}.
\end{equation}

Now, let us introduce new, rescaled coordinates $(X_1,Y_1,X_2,Y_2)$
and parameters $(B,E_0,\dots,E_{m-1})$ by the following rule:
\begin{equation}\label{pochti}
\begin{array}{l} \displaystyle
x_1 = b_2 d_1 \gamma_1^{-k_1/m}\gamma_2^{-k_2/m} X_1,\qquad
x_2=D b_2 \lambda_1^{k_1}\gamma_1^{-k_1/m}\gamma_2^{-k_2/m} X_2,\\ \displaystyle
y_1=\gamma_1^{-k_1/m}\gamma_2^{-k_2(1+1/m)} Y_1, \qquad
y_2=d_1 \gamma_1^{-k_1/m}\gamma_2^{-k_2/m} Y_2,\\
\displaystyle
\tilde\mu_s = d_1^s \left(\gamma_1^{k_1/m}\gamma_2^{k_2/m)}\right)^{-(m+1-s)} E_s \;\;\;\;
(s=0,\dots,m-1),\\ \displaystyle
\theta=\beta_0+\frac{\ln(-B/(Db_2c_2)}{k_2|\ln|\lambda_1\gamma_1||};
\end{array}
\end{equation}
recall that $Db_2c_2<0$ by our assumptions, so the new parameter $B$ should be positive. After
the rescaling, the maps $T_{12}T_{01}^{k_1}$ (given by (\ref{tlm})) and $T_{21}T_{02}^{k_2}$
(given by (\ref{tlm2})) are rewritten as follows (we take into account
(\ref{kb}),(\ref{lgin})):
$$\begin{array}{l}\displaystyle
\bar X_2 = X_1 + o(1)_{k_{1,2}\rightarrow+\infty}\\ \displaystyle
Y_1= \bar Y_2 + o(1)_{k_{1,2}\rightarrow+\infty}, \quad \mbox{ and}\\
\bar X_1 = Y_2 + o(1)_{k_{1,2}\rightarrow+\infty},\\
\bar Y_1= -B X_2 +
\sum_{s=0}^{m-1} E_s Y_2^s + d\;Y_2^{m+1}+o(1)_{k_{1,2}\rightarrow+\infty},
\end{array}
$$
where $d:=d_2 d_1^{m+1}\neq0$. It follows
that with our choice of $k_{1,2}$ the map ${\cal T}:=T_{21}T_{02}^{k_2}T_{12}T_{01}^{k_1}$
on $\sigma_{1k_1}$ can be written as
$$\bar X = Y + o(1)_{k_{1,2}\rightarrow+\infty},\qquad
\bar Y= -B X +
\sum_{s=0}^{m-1} E_s Y^s + d\;Y^{m+1}+o(1)_{k_{1,2}\rightarrow+\infty},
$$ where $(X,Y)$ are the coordinates which were previously denoted as $(X_1,Y_1)$. By denoting
the right-hand side of the equation for $\bar X$ as the new $Y$-variable, we finally write
the map $\cal T$ in the following form
\begin{equation}\label{retmap}
\bar X = Y,\qquad
\bar Y= -B X +
\sum_{s=0}^{m-1} E_s Y^s + d\;Y^{m+1}+o(1)_{k_{1,2}\rightarrow+\infty}.
\end{equation}
The transformation to
coordinates $(X,Y)$ was of class $C^\rho$ with respect to the original coordinates and $C^{\rho-m}$
with respect to the parameters (more precisely: the $m$-th derivatives with respect to the coordinates
are $C^{\rho-m}$ with respect to both the coordinates and parameters). The $o(1)$-terms in (\ref{retmap})
are functions of $(X,Y)$ and $(B,E_0,\dots,E_{m-1})$ which tend to zero as $k_{1,2}\rightarrow+\infty$
along such derivatives up to the order $\rho$, for which the number of differentiations with respect
to the parameters does not exceed $(\rho-m)$. In any case, since $\rho\geq m+1$,
we have that their derivatives with respect to $(X,Y)$ up to the order $m$ tend to zero, all with at
least one derivative with respect to $(B,E_0,\dots,E_{m-1})$.

We mention that a similar form for the rescaled first-return map was obtained in \cite{St} for the
case $m=1$. Map (\ref{retmap}) was also obtained in \cite{GTS6} as the rescaled first-return map
near a homoclinic tangency of order $m$. Let us show
that if the $o(1)$-terms are sufficiently small, then a map of form (\ref{retmap})
has, at certain uniformly bounded values of $E$ and $B>0$ and in a bounded region of $(X,Y)$,
a fixed point near which the map is smoothly conjugate
to (\ref{flat}). This will give us the lemma: by (\ref{pochti}),(\ref{shifto}),
any bounded values of $(X,Y)$
correspond to $(x_1,y_1)\rightarrow(x_1^+,y_1^-)$ as $k_{1,2}\rightarrow+\infty$, i.e. if the values
of $(X,Y)$ are
bounded at the fixed point, then the latter belongs
to the domain of definition of the first-return map $T_{21}T_{02}^{k_2}T_{12}T_{01}^{k_1}$
at all sufficiently large $k_{1,2}$, and any bounded values of $(E_0,\dots,E_{m-1})$ and $B>0$
correspond to $\nu_k\rightarrow0$ as $k_{1,2}\rightarrow+\infty$.

We may take any $Y^*$ as the coordinate $X=Y=Y^*$ of the fixed point of (\ref{retmap}); then
the second equation defines the corresponding value of $E_0=Y^*(1+B)-
\sum_{s=1}^{m-1} E_s (Y^*)^s - d\;(Y^*)^{m+1}+o(1)$. After shifting the origin to the fixed point,
the map takes the form
\begin{equation}\label{rets}
\begin{array}{l}
\bar X = Y,\\
\bar Y= -\hat B X + 2Y+
\sum_{s=1}^{m} \hat E_s Y^s + O(Y^{m+1})+O(|XY|+X^2)\cdot o(1)_{k\rightarrow+\infty},
\end{array}
\end{equation}
where the coefficients $\hat E_s$, $s=1,\dots,m-1$, are related to $E_1,\dots,E_{m-1}$ by a
diffeomorphism:
$$\hat E_1=E_1-2+\sum_{j=2}^{m-1} j E_j (Y^*)^{j-1}+(m+1)d\;(Y^*)^m+o(1),$$
$$\hat E_s=E_s+\sum_{j=s+1}^{m-1}E_j C^s_j (Y^*)^{j-s}+C^s_{m+1}d\;(Y^*)^{m+1-s}+o(1) \; (s=2,\dots,m-1),$$
and
$$\hat E_m=(m+1)d\;Y^*+o(1),\qquad \hat B=B+o(1).$$
As we see, bounded values of $\hat B,\hat E_1,\dots,\hat E_m$
give us bounded values of $B,E_0,\dots,E_{m-1}$ as well.

We will further fix $\hat B=1$ and will take $(\hat E_1,\dots,\hat E_m)=o(1)_{k\rightarrow+\infty}$, so
the map will limit, as $k\rightarrow+\infty$, to
$$\bar X=Y, \qquad \bar Y =-X+2Y.$$
The fixed point of this map at $(X,Y)=0$ has a double multiplier $(+1)$. Let us recall the normal
form theory for maps with such a fixed point. In general, we can always write such map in the form
\begin{equation}\label{parnorf}
\begin{array}{l}\displaystyle
\bar u=u+v+\sum_{2\leq i+j\leq m} \alpha_{ij} u^iv^j +o(|u|^m+|v|^m),\\ \displaystyle
\bar v = v + \sum_{2\leq i+j\leq m} \beta_{ij} u^iv^j +o(|u|^m+|v|^m).
\end{array}
\end{equation}
By denoting $z=v+\sum_{2\leq i+j\leq m} \alpha_{ij} u^iv^j +o(|u|^m+|v|^m)\equiv \bar u-u$,
the map takes the form
\begin{equation}\label{parnorf2}
\bar u=u+z,
\qquad \bar z = z + \sum_{2\leq i+j\leq m} \phi_{ij} u^iz^j +o(|u|^m+|z|^m).
\end{equation}
Here
\begin{equation}\label{norfgam}
\phi_{ij}=\beta_{ij}+\sum_{i<k\leq i+j} C_k^i\alpha_{k,i+j-k}+G_{ij},
\end{equation}
where $G_{ij}$ stand for polynomial functions of the
coefficients $\alpha_{i^{\prime}j^{\prime}}$ and $\beta_{i^{\prime}j^{\prime}}$ with
$i^{\prime}+j^{\prime}<i+j$; each monomial in $G_{ij}$ has degree $2$ or higher.
After the polynomial coordinate transformation $(u,z)\mapsto(U,Z)$, where
$$U=u+\sum_{0\leq j\leq n-2} A_j u^{n-j} z^j,\qquad Z=z+
\sum_{0\leq j\leq n-2} A_j (\bar u^{n-j} \bar z^j-u^{n-j} z^j),$$
the map will retain its form (\ref{parnorf2}) with the coefficients $\phi_{ij}$ unchanged
for $i+j<n$, and
$\phi_{ij,new}=\phi_{ij}+\sum_{0\leq l\leq j-2} (2^{j-l}-2) C^{j-l}_{n-l} A_l$
at $i+j=n$. As we see, by choosing $A_l$ appropriately, one can make $\phi_{ij}$
with $j\geq2, i=n-j$ assume any given values. Thus, by making such transformations consecutively,
with $n$ running from $2$ to $m$, one can make all the coefficients $\phi_{ij}$ with $j\geq 2$
vanish; the new coefficients $\phi_{ij}$ ($j=0,1$) differ from the old ones by nonlinear (i.e.
quadratic and higher order) expressions
which depend only on the coefficients $\phi_{i^{\prime}j^{\prime}}$ with $i^{\prime}+j^{\prime}<i+j$,
i.e. they are still given by (\ref{norfgam}) with some new $G_{ij}$.
Obviously, for any fixed $m$, we arrive to a similar conclusion for any map (with a fixed point at zero)
whose linear part is sufficiently close to that of (\ref{parnorf}). Specifically,
for any sufficiently small $\varepsilon_{1,2}$, a map of the form
\begin{equation}\label{parnf}
\begin{array}{l}\displaystyle
\bar u=u+v+\sum_{2\leq i+j\leq m} \alpha_{ij} u^iv^j +o(|u|^m+|v|^m),\\ \displaystyle
\bar v = (1+\varepsilon_1)v +\varepsilon_2 u + \sum_{2\leq i+j\leq m} \beta_{ij} u^iv^j +o(|u|^m+|v|^m),
\end{array}
\end{equation}
can be brought by a smooth coordinate transformation to the normal form
\begin{equation}\label{parnf2}
\bar u=u+z,
\qquad \bar z = (1+\varepsilon_1)z + \varepsilon_2 u + \sum_{i=2}^m (\phi_{i0} u^i+\phi_{i-1,1}u^iz)
+o(|u|^m+|z|^m),
\end{equation}
where
\begin{equation}\label{norfpam}
\phi_{i0}=\beta_{i0}+G_{i0}+H_{i0},\qquad
\phi_{i1}=\beta_{i1}+(i+1)\alpha_{i+1,0}+G_{i1}+H_{i1}.
\end{equation}
Here $G_{ij}$ is a polynomial of $\alpha_{i^{\prime}j^{\prime}}$ and
$\beta_{i^{\prime}j^{\prime}}$ with  $i^{\prime}+j^{\prime}<i+j$, all monomials of $G_{ij}$ are
of degree at least $2$;
the function $H_{ij}$ is a linear function of $\alpha_{i^{\prime}j^{\prime}}$ and
$\beta_{i^{\prime}j^{\prime}}$ with  $i^{\prime}+j^{\prime}=i+j$; the coefficients of $G_{ij}$ and
$H_{ij}$ depend on $\varepsilon_{1,2}$ as well, and $H_{ij}$ vanish at $\varepsilon_1=\varepsilon_2=0$.

In fact, given any $m$, for any sufficiently small $\varepsilon_{1,2}$, there exists a pair of functions
$\psi_{0,1}(u)$ such that
after an appropriate choice of the coordinates $(u,v)$, map (\ref{parnf}) can
be made $o(|u|^m+|v|^m)$-close to the time-$1$ map by a flow of the form
\begin{equation}\label{flwuv}
\dot u=v, \qquad \dot v= \psi_0(u)+v\psi_1(u).
\end{equation}
Obviously, the flow has to have an equilibrium at zero and the eigenvalues of the corresponding
linearization matrix have to be logarithms of the eigenvalues of
$\left(\begin{array}{cc} 1 & 1\\ \varepsilon_2 & 1+\varepsilon_1\end{array}\right)$, the linearization
matrix of (\ref{parnf}). Thus, we assume
\begin{equation}\label{psilin12}
\begin{array}{l}\displaystyle
\psi_0(0)=0,\qquad \psi_1(0)=\ln(1+\varepsilon_1-\varepsilon_2),\\\displaystyle
\psi_0'(0)=-\ln(1+\frac{\varepsilon_1}{2}+\sqrt{\varepsilon_2+\frac{\varepsilon_1^2}{4}})
\ln(1+\frac{\varepsilon_1}{2}-\sqrt{\varepsilon_2+\frac{\varepsilon_1^2}{4}}).
\end{array}
\end{equation}
The time-$t$ map of (\ref{flwuv}) can be found by expanding it in powers of the initial conditions
$(u(0),v(0))$. This
gives us the time-$1$ map in the form which is brought, by an $O(\varepsilon_{1,2})$-close to identity
linear transformation, to form (\ref{parnf}) with
\begin{equation}\label{flalb}
\begin{array}{l}\displaystyle
\alpha_{ij}=\frac{1}{j+1}\left(\psi_{i+j,0}\left[\frac{C^{i}_{i+j}}{j+2}+O(\varepsilon_{1,2})\right]+
\psi_{i+j-1,1}\left[\frac{C^{i}_{i+j-1}}{j}+O(\varepsilon_{1,2})\right]\right)+\tilde\alpha_{ij},\\\displaystyle
\beta_{ij}=\psi_{i+j,0}\left[\frac{C^{i}_{i+j}}{j+1}+O(\varepsilon_{1,2})\right]+
\psi_{i+j-1,1}\left[\frac{C^{i}_{i+j-1}}{j}+O(\varepsilon_{1,2})\right]+\tilde\beta_{ij},
\end{array}
\end{equation}
where $\tilde\alpha_{ij}$ and $\tilde\beta_{ij}$ are nonlinear functions of the coefficients $\psi_{i'j'}$
with $i'+j'<i+j$ (we denote $\psi_j(u)=\sum_i\psi_{ij}u^i$). Two maps of form (\ref{parnf})
can be made $o(|u|^m+|v|^m)$-close by means of a smooth coordinate transformation if their
normal forms (\ref{parnf2}) coincide up to the order $m$. Thus, we find from formulas
(\ref{norfpam}),(\ref{psilin12}),(\ref{flalb}) that a map of form (\ref{parnf}) can indeed, after
an appropriate coordinate transformation, be made $o(|u|^m+|v|^m)$-close
to the time-$1$ map by a flow of form (\ref{flwuv}) with the coefficients defined by
\begin{equation}\label{flps}
\psi_{i0}=\beta_{i0}+\tilde\psi_{i0},\qquad
\psi_{i1}=\beta_{i1}+(i+1)(\alpha_{i+1,0}-\beta_{i+1,0})+\tilde\psi_{i1},
\end{equation}
where $\tilde\psi_{ij}$ is a polynomial function of $\alpha_{i'j'}$
and $\beta_{i'j'}$ with $i'+j'\leq i+j$; the monomials that include $\alpha_{i'j'}$
or $\beta_{i'j'}$ with $i'+j'<i+j$ are of degree $2$ or higher, while the terms that
include $\alpha_{i'j'}$ or $\beta_{i'j'}$ with $i'+j'=i+j$ are linear and the corresponding
coefficients are small of order $O(\varepsilon_{1,2})$. It follows from
(\ref{psilin12}) that formulas (\ref{flps}) remain valid at $j=1$,
if we formally put $\beta_{10}=\varepsilon_2$, $\beta_{01}=\varepsilon_1$, $\alpha_{10}=\alpha_{01}=0$.

Returning to map (\ref{rets}), by putting $X=u$ and $Y=u+v$, we bring the map
at $\hat B=1$ to form
$$\bar u = u+v,\qquad\bar v= v + \beta_0(u)+v\beta_1(u)+O(v^2)+o(|u|^m+|v|^m),$$
where the coefficients of the polynomials $\beta_{j}=\sum_{i=1}^{m-j}\beta_{ij}u^i$ ($j=0,1$)
satisfy
$$ \beta_{i0}=\hat E_i+o(1)_{k\rightarrow+\infty}, \qquad
\beta_{i1}=(i+1)\hat E_{i+1}+o(1)_{k\rightarrow+\infty}.
$$
Thus, by making normal form transformations described
above, we can, at $\hat E_1=\varepsilon$ sufficiently small, make the map $o(|u|^m+|v|^m)$-close to
the time-$1$ map of the flow (\ref{flwuv}) with
$$\psi_{i0}=\hat E_i+o(1)_{k\rightarrow+\infty}+o(\hat E_1,\dots,\hat E_m),\qquad
\psi_{i1}=o(1)_{k\rightarrow+\infty}+o(\hat E_1,\dots,\hat E_m);$$
see (\ref{flps}). It follows that at sufficiently large $k$ we can always choose the values
of $\hat E_1,\dots,\hat E_m$ in such a way that $\psi_0(u)\equiv s\psi_1(u)$, where
$s=1$ if $\psi_1'(0)\leq0$, and $s=-1$ if $\psi_1'(0)>0$.

As we see, at an appropriate choice of the parameters the map can be made
$o(|u|^m+|v|^m)$-close to the time-$1$ map of the flow
$$\dot u=v,\qquad \dot v = \psi_1(u)(s+v),$$
which takes the desired form (\ref{flowat}) after the change $(u,v)\rightarrow(su,sv)$; here
$\Psi(u)=-\psi_1(su)$, so the requirement $\Psi'(0)\geq 0$ is ensured by our choice of $s$.~$\Box$\\

{\small This work was supported by grant ISF 273/07 and by Royal Society grant
``Homoclinic bifurcations'', and it was done during my stay at MPIM, Bonn.}

\end{document}